\documentclass{biometrika}

\usepackage{times}
\usepackage{bm}
\usepackage{natbib}

\usepackage{graphicx}
\usepackage{amssymb}
\usepackage[cmex10]{amsmath}

\begin{document}

\received{November 2010}
\revised{April 2011}

\markboth{D. S. Choi, P. J. Wolfe \and E. M. Airoldi}{Stochastic blockmodels}

\title{Stochastic blockmodels with growing number of classes}

\author{D. S. CHOI}
\affil{School of Engineering and Applied Sciences, Harvard University, Cambridge, Massachusetts 02138, U.S.A. \email{dchoi@seas.harvard.edu}}
\author{P. J. WOLFE}
\affil{School of Engineering and Applied Sciences, and Department of Statistics, Harvard University, Cambridge, Massachusetts 02138, U.S.A. \email{wolfe@stat.harvard.edu}}
\author{\and E. M. AIROLDI}
\affil{Department of Statistics, Harvard University, Cambridge, Massachusetts 02138, U.S.A. \email{airoldi@fas.harvard.edu}}

\maketitle

\begin{abstract}
We present asymptotic and finite-sample results on the use of stochastic blockmodels for the analysis of network data.  We show that the fraction of misclassified network nodes converges in probability to zero
under maximum likelihood fitting when the number of classes is allowed to grow as the root of the network size and the average network degree grows at least poly-logarithmically in this size.  We also establish finite-sample confidence bounds on maximum-likelihood blockmodel parameter estimates from data comprising independent Bernoulli random variates; these results hold uniformly over class assignment.  We provide simulations verifying the conditions sufficient for our results, and conclude by fitting a logit parameterization of a stochastic blockmodel with covariates to a network data example comprising a collection of Facebook profiles, resulting in block estimates that reveal residual structure.
\end{abstract}


\section{Introduction}

The global structure of social, biological, and information networks is sometimes envisioned as the aggregate of many local interactions whose effects propagate in ways that are not yet well understood. There is increasing opportunity to collect data on an appropriate scale for such systems, but their analysis remains challenging \citep{goldenberg2009survey}. Here we analyze a statistical model for network data known as the (single-membership) stochastic blockmodel. Its salient feature is that it partitions the $N$ nodes of a network into $K$ distinct classes whose members all interact similarly with the network.  Blockmodels were first associated with the deterministic concept of structural equivalence in social network analysis \citep{lorrain1971structural}, where two nodes were considered interchangeable if their connections were equivalent in a formal sense. This concept was adapted to stochastic settings and gave rise to the stochastic blockmodel in work by \citet{holland1983stochastic} and \citet{fienberg1985statistical}.  The model and extensions thereof have since been applied in a variety of disciplines \citep{wang1987stochastic, nowicki2001estimation, girvan2002community, airoldi2005latent, doreian2005generalized, newman2006modularity, handcock2007model, hoff2008modeling, airoldi2008mixed, copic2009identifying, mariadassou2010uncovering, karrer2011stochastic}.

In this work we provide a finite-sample confidence bound that can be used when estimating network structure from data modeled by independent Bernoulli random variates, and also show that under maximum likelihood fitting of a correctly specified $K$-class blockmodel, the fraction of misclassified network nodes converges in probability to zero even when the number of classes $K$ grows with $N$.  As noted by \citet{rohe2010spectral}, this is advantageous if we expect class sizes to remain relatively constant even as $N$ increases.  Related results for fixed $K$ have been shown by \citet{snijders1997estimation} for networks with linearly increasing degree, and in a stronger sense for sparse graphs with poly-logarithmically increasing degree by \citet{bickel2009nonparametric}.

Our results can be related to those of \citet{rohe2010spectral}, who use spectral methods to bound the number of misclassified nodes in the stochastic blockmodel with increasing $K$, although with the more restrictive requirement of nearly linearly increasing degree.  As noted by those authors, this assumption may not hold in many practical settings.  Our manner of proof requires only poly-logarithmically increasing degree, and is more closely related to the fixed-$K$ proof of \citet{bickel2009nonparametric}, although we note that spectral clustering as suggested by \citet{rohe2010spectral} provides a computationally appealing alternative to maximum likelihood fitting in practice.

As discussed by \citet{bickel2009nonparametric}, one may assume exchangeability in lieu of a generative $K$-class blockmodel: An analogue to de~Finetti's theorem for exchangeable sequences states that the probability distribution of an infinite exchangeable random graph is expressible as a mixture of distributions whose components can be approximated by blockmodels \citep{kallenberg2005probabilistic, bickel2009nonparametric}.  An observed network can then be viewed as a sample drawn from this infinite conceptual population, and so in this case the fitted blockmodel describes one mixture component thereof.

\section{Statement of results}

\subsection{Problem formulation and definitions}

We consider likelihood-based inference for independent Bernoulli data $\{A_{ij}\}\ (i = 1, \ldots, N; j = i + 1, \ldots, N)$, both when no structure linking the success probabilities $\{P_{ij}\}$ is assumed, as well as the special case when a stochastic blockmodel of known order $K$ is assumed to apply.  To this end, let $A \in \{0,1\}^{N \times N}$ denote the symmetric adjacency matrix of a simple, undirected graph on $N$ nodes whose entries $\{A_{ij}\}$ for $i < j$ are assumed independent $\operatorname{Bernoulli}(P_{ij})$ random variates, and whose main diagonal $\{A_{ii}\}_{i=1}^N$ is fixed to zero.  The average degree of this graph is $2M/N$, where $M =  \sum_{i<j} P_{ij}$ is its expected number of edges.  Under a $K$-class stochastic blockmodel, these edge probabilities are further restricted to satisfy
\begin{equation}
\label{eq:Pij}
P_{ij} = \theta_{z_i z_j} \quad (i = 1, \ldots, N; j = i + 1, \ldots, N)
\end{equation}
for some symmetric matrix $\theta \in [0,1]^{K\times K}$ and membership vector $z \in \{1,\ldots,K\}^N$.  Thus the probability of an edge between two nodes is assumed to depend only on the class of each node.

Let $L(A;z,\theta)$ denote the log-likelihood of observing data matrix $A$ under a $K$-class blockmodel with parameters $(z,\theta)$, and $\bar{L}_P(z,\theta)$ its expectation:
\begin{align*}
    L(A;z,\theta) & = \sum_{i<j} \left\{ A_{ij} \log \theta_{z_iz_j} + (1-A_{ij})\log(1-\theta_{z_iz_j}) \right\}, \\
\bar{L}_P(z,\theta) & = \sum_{i<j} \left\{ P_{ij} \log \theta_{z_iz_j} + (1-P_{ij})\log(1-\theta_{z_iz_j}) \right\}.
\end{align*}

For fixed class assignment $z$, let $N_a$ denote the number of nodes assigned to class $a$, and let $n_{ab}$ denote the maximum number of possible edges between classes $a$ and $b$; i.e., $n_{ab} = N_a N_b$ if $a \neq b$ and $n_{aa} = {N_a \choose 2}$.  Further, let $\hat{\theta}^{(z)}$ and $\bar{\theta}^{(z)}$ be symmetric matrices in $[0,1]^{K \times K}$, with
\begin{align*}
\hat{\theta}_{ab}^{(z)} &= \frac{1}{n_{ab}}\sum_{i<j} A_{ij} \, 1\{z_i=a,z_j=b\}
\quad (a = 1, \ldots, K; b = a, \ldots, K), \\
\bar{\theta}_{ab}^{(z)} &= \frac{1}{n_{ab}}\sum_{i<j} P_{ij} \, 1\{z_i=a,z_j=b\}
\quad (a = 1, \ldots, K; b = a, \ldots, K)
\end{align*}
defined whenever $n_{ab} \neq 0$.  Observe that $\hat{\theta}^{(z)}$ comprises sample proportion estimators as a function of $z$, whereas $\bar{\theta}^{(z)}$ is its expectation under the independent $\{\operatorname{Bernoulli}(P_{ij})\}$ model.  Taken over all class assignments $z \in \{1,\ldots,K\}^N$, the sets $\{\hat{\theta}^{(z)}\}$ comprise a sufficient statistic for the family of $K$-class stochastic blockmodels, and for each $z$, $\hat{\theta}^{(z)}$ maximizes $L(A;z,\cdot)$. Analogously, the sets $\{\bar{\theta}^{(z)}\}$ are functions of the model parameters $\{P_{ij}\}_{i<j}$, and maximize $\bar{L}_P(z,\cdot)$.  We write $\hat{\theta}$ and $\bar{\theta}$ when the choice of $z$ is understood, and $L(A;z)$ and $\bar{L}_P(z)$ to abbreviate $\sup_\theta L(A;z,\theta)$ and $\sup_\theta \bar{L}_P(z,\theta)$ respectively.

Finally, observe that when a blockmodel with parameters $(\bar{z}, \bar{\theta})$ is in force, then $P_{ij} = \bar{\theta}_{\bar{z}_i \bar{z}_j}$ in accordance with~\eqref{eq:Pij}, and consequently $\bar{L}_P$ is maximized by the true parameter values $(\bar{z}, \bar{\theta})$:
\begin{equation*}
\bar{L}_P(\bar{z}, \bar{\theta}) - \bar{L}_P(z,\theta) = \sum_{i<j} D(P_{ij} \mid\mid \theta_{z_i z_j } ) \geq \sum_{i<j} 2(P_{ij}-\theta_{z_i z_j})^2 \geq 0,
\end{equation*}
where $D(p \mid\mid p')$ denotes the Kullback--Leibler divergence of a $\operatorname{Bernoulli}(p')$ distribution from a $\operatorname{Bernoulli}(p)$ one.

\subsection{Fitting a $K$-class stochastic blockmodel to independent Bernoulli trials}

Fitting a $K$-class stochastic blockmodel to independent $\operatorname{Bernoulli}(P_{ij})$ trials yields estimates $\hat{\theta}^{(z)}$ of averages $\bar{\theta}^{(z)}$ of subsets of the parameter set $\{P_{ij}\}$, with each class assignment $z$ inducing a partition of that set.  We begin with a basic lemma that expresses the difference $L(A;z) - \bar{L}_P(z)$ in terms of $\hat{\theta}^{(z)}$ and $\bar{\theta}^{(z)}$, and follows directly from their respective maximizing properties.

\begin{lemma}\label{lem:LbarLDiff}
Let $\{A_{ij}\}_{i< j}$ comprise independent $\operatorname{Bernoulli}(P_{ij})$ trials.  Then the difference $\sup_\theta L(A;z,\theta) - \sup_\theta \bar{L}_P(z,\theta)$ can be expressed for $X = \sum_{i<j} A_{ij} \log \{\bar{\theta}_{z_i z_j} / (1-\bar{\theta}_{z_iz_j})\}$ as
\begin{equation*}
L(A;z) - \bar{L}_P(z) = \textstyle \sum_{a\leq b} n_{ab} D(\hat{\theta}_{ab} \mid\mid \bar{\theta}_{ab}) + X - E(X).
\end{equation*}
\end{lemma}

We first bound the former quantity in this expression, which provides a measure of the distance between $\hat{\theta}$ and its estimand $\bar{\theta}$ under the setting of Lemma~\ref{lem:LbarLDiff}.  The bound is used in subsequent asymptotic results, and also yields a kind of confidence measure on $\hat{\theta}$ in the finite-sample regime.

\begin{theorem}\label{thm:sumKLDivBnd}
Suppose that a $K$-class stochastic blockmodel is fitted to data $\{A_{ij}\}_{i< j}$ comprising $\binom{N}{2}$ independent $\operatorname{Bernoulli}(P_{ij})$ trials, where, for any class assignment $z$, estimate $\hat{\theta}$ maximizes the blockmodel log-likelihood $L(A;z,\cdot)$.  Then with probability at least $1-\delta$,
\begin{equation}\label{eq:sumKLDivBnd}
\max_z \left\{ \textstyle \sum_{a\leq b} n_{ab} D(\hat{\theta}_{ab} \mid\mid \bar{\theta}_{ab}) \right \} < N \log K + (K^2+K)\log\Big(\frac{N}{K}+1\Big) + \log\frac{1}{\delta}.
\end{equation}
\end{theorem}

Theorem~\ref{thm:sumKLDivBnd} is proved in the Appendix via the method of types: for fixed $z$, the probability of any realization of $\hat{\theta}$ is first bounded by $\exp\{-\sum_{a\leq b} n_{ab} D(\hat{\theta}_{ab} \mid\mid \bar{\theta}_{ab})\}$.  A counting argument then yields a deviation result in terms of $(N/K+1)^{K^2+K}$, and finally a union bound is applied so that the result holds uniformly over all $K^N$ possible choices of assignment vector $z$.

Our second result is asymptotic, and combines Theorem~\ref{thm:sumKLDivBnd} with a Bernstein inequality for bounded random variables, applied to the latter terms $X - E(X)$ in Lemma~\ref{lem:LbarLDiff}.  To ensure boundedness we assume minimal restrictions on each $P_{ij}$; this Bernstein inequality, coupled with a union bound to ensure that the result holds uniformly over all $z$, dictates growth restrictions on $K$ and $M$.

\begin{theorem}\label{thm:unifLikBnd}
Assume the setting of Theorem~\ref{thm:sumKLDivBnd}, whereby a $K$-class blockmodel is fitted to $\binom{N}{2}$ independent $\operatorname{Bernoulli}(P_{ij})$ random variates $\{A_{ij}\}_{i< j}$, and further assume that $1/N^2 \leq P_{ij} \leq 1 - 1/N^2$ for all $N$ and $i<j$.  Then if $K = \mathcal{O}(N^{1/2})$ and $M = \omega(N (\log N)^{3+\delta})$ for some $\delta > 0$,
$$
\max_z |L(A;z) - \bar{L}_P(z)| = o_P(M).
$$
\end{theorem}
Thus whenever each $P_{ij}$ is bounded away from 0 and 1 in the manner above, the maximized log-likelihood function $L(A;z) = \sup_\theta L(A;z,\theta)$ is asymptotically well behaved in network size $N$ as long as the network's average degree $2M/N$ grows faster than $(\log N)^{3+\delta}$ and the number $K$ of classes fitted to it grows no faster than $N^{1/2}$.

\subsection{Fitting a correctly specified $K$-class stochastic blockmodel}

The above results apply to the general case of independent Bernoulli data $\{A_{ij}\}$, with no additional structure assumed amongst the set of success probabilities $\{P_{ij}\}$; if we further assume the data to be generated by a $K$-class stochastic blockmodel whose parameters $(\bar{z},\bar{\theta})$ are subject to suitable identifiability conditions, it is possible to characterize the behavior of the class assignment estimator $\hat{z}$ under maximum likelihood fitting of a correctly specified $K$-class blockmodel.

\begin{theorem}\label{thm:classConv}
If the conclusion $\max_z |L(A;z) - \bar{L}_P(z)| = o_P(M)$ of Theorem~\ref{thm:unifLikBnd} holds, and data are generated according to a $K$-class blockmodel with membership vector $\bar{z}$, then
\begin{equation}\label{eq:likeliConv}
\bar{L}_P(\bar{z}) - \bar{L}_P(\hat{z}) = o_P(M),
\end{equation}
with respect to the maximum-likelihood $K$-class blockmodel class assignment estimator $\hat{z}$.

Let $N_{\mathrm{e}}(\hat{z})$ be the number of incorrect class assignments under $\hat{z}$, counted for every node whose true class under $\bar{z}$ is not in the majority within its estimated class under $\hat{z}$.  If furthermore the following identifiability conditions hold with respect to the model sequence:

\emph{(i)} for all blockmodel classes $a = 1,\ldots,K$, class size $N_a$ grows as $\min_a \{N_a\} = \Omega(N/K)$;

\emph{(ii)} the following holds over all distinct class pairs $(a,b)$ and all classes $c$:
$$
\min_{(a,b)}\, \max_{c} \Big\{
D\Big(\bar{\theta}_{ac} \mid\mid \frac{\bar{\theta}_{ac}+\bar{\theta}_{bc}}{2} \Big) +
D\Big(\bar{\theta}_{bc} \mid\mid \frac{\bar{\theta}_{ac}+\bar{\theta}_{bc}}{2} \Big)
\Big \}
= \Omega\Big(\frac{MK}{N^2}\Big),
$$
then it follows from~\eqref{eq:likeliConv} that
$
N_{\mathrm{e}}(\hat{z}) = o_P(N).
$
\end{theorem}

Thus the conclusion of Theorem~\ref{thm:classConv} is that under suitable conditions the fraction $N_{\mathrm{e}}/N$ of misclassified nodes goes to zero in $N$, yielding a convergence result for stochastic blockmodels with growing number of classes.  Condition (i) stipulates that all class sizes grow a rate that is eventually bounded below by a single constant times $N/K$, while condition (ii) ensures that any two rows of $\theta$ differ in at least one entry by an amount that is eventually bounded by a single constant times $MK/N^2$.  Observe that if eventually $K = N^{1/2}$ and $M = N (\log N)^4$ so that conditions on $K$ and $M$ sufficient for Theorem~\ref{thm:unifLikBnd} are met, then since $(\log N)^4 = o(N^{1/2})$, it follows that $MK/N^2$ goes to zero in $N$.

\section{Numerical results}

We now present results of a small simulation study undertaken to investigate the assumptions and conditions of Theorems~\ref{thm:sumKLDivBnd}--\ref{thm:classConv} above, in which $K$-class blockmodels were fitted to various networks generated at random from models corresponding to each of the three theorems.  Because exact maximization in $z$ of the blockmodel log-likelihood $L(A;z,\theta)$ is computationally intractable even for moderate $N$, we instead employed Gibbs sampling to explore the function $\max_\theta L(A;z,\theta)$ and recorded the best value of $z$ visited by the sampler.  As the results of Theorems~\ref{thm:sumKLDivBnd} and~\ref{thm:unifLikBnd} hold uniformly in $z$, however, we expect $\bar{\theta}$ and $\bar{L}_P(z)$ to be close to their empirical estimates whenever $N$ is sufficiently large, regardless of the approach employed to select $z$.  This fact also suggests that a single-class (Erd\"{o}s-R\'{e}nyi) blockmodel may come closest to achieving equality in Theorems~\ref{thm:sumKLDivBnd} and~\ref{thm:unifLikBnd}, as many class assignments are equally likely a priori to have high likelihood.  By similar reasoning, a weakly identifiable model should come closest to achieving the error bound in Theorem~\ref{thm:classConv}, such as one with nearly identical within- and between-class edge probabilities.  We describe each of these cases empirically in the remainder of this section.

First, the tightness of the confidence bound of~\eqref{eq:sumKLDivBnd} from Theorem~\ref{thm:sumKLDivBnd} was investigated by fitting $K$-class blockmodels to Erd\"{o}s-R\'{e}nyi networks comprising $\binom{N}{2}$ independent $\operatorname{Bernoulli}(p)$ trials, with $N = 500$ nodes and $p =\ $0$\cdot$075 chosen to match the data analysis example in the sequel, and $K \in \{5,10,20,30,40,50\}$.
For each $K$, the error terms $\sum_{a\leq b} n_{ab} D(\hat{\theta}_{ab} \mid\mid \bar{\theta}_{ab})$ and $\{\sum_{a \leq b} n_{ab} (\hat{\theta}_{ab} - \bar{\theta}_{ab})^2\}^{1/2}$ were recorded for each of 100 trials and compared to the respective 95\% confidence bounds ($\delta =\ $0$\cdot$05) derived from Theorem~\ref{thm:sumKLDivBnd}.  The bounds overestimated the respective errors by a factor of 3 to 7 on average, with small standard deviation.  In this worst-case scenario the bound is loose, but not unusable; the errors never exceeded the 95\% confidence bounds in any of the trials.

\begin{figure}
\begin{center}
\includegraphics[width=.34\columnwidth]{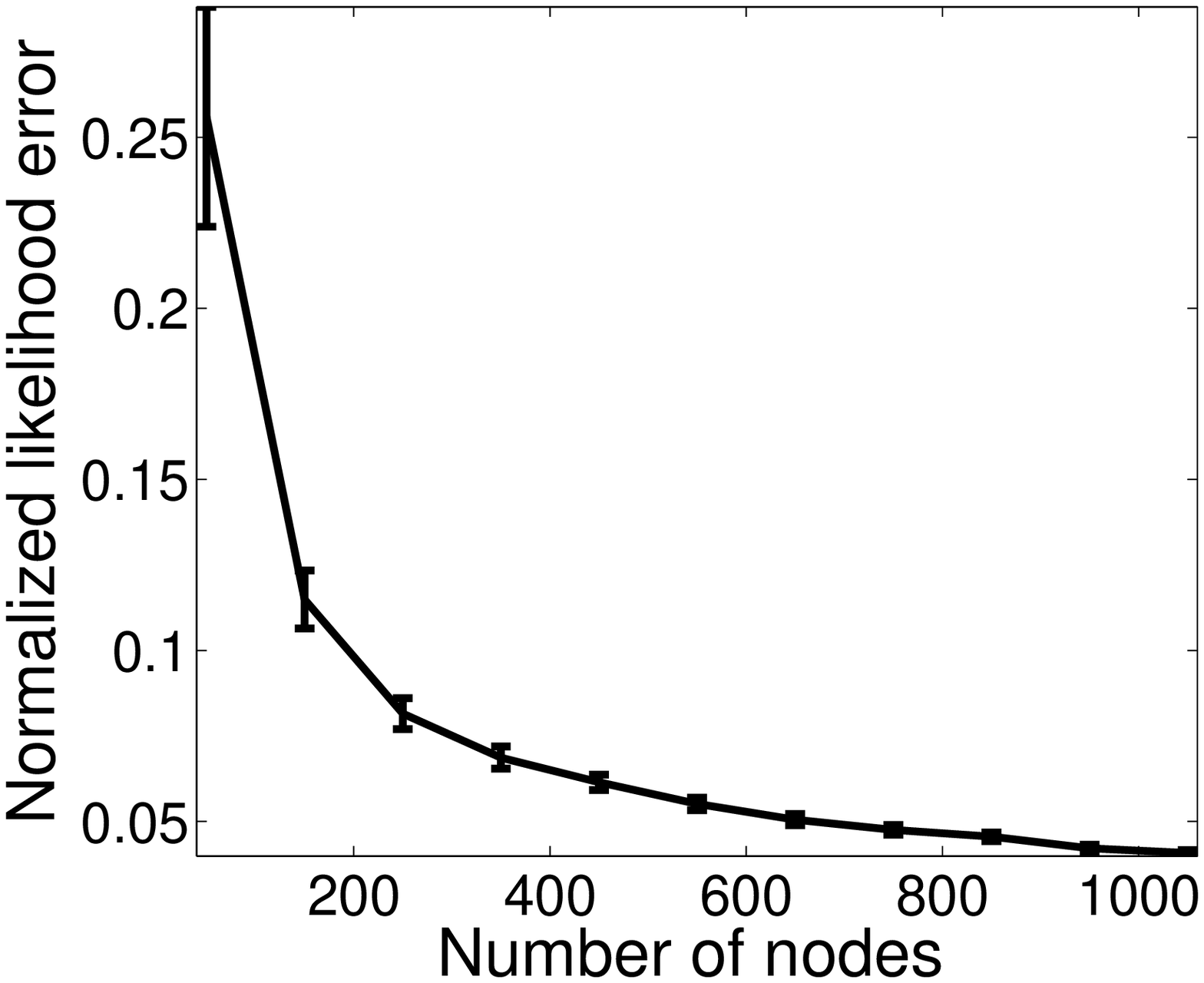}\hspace{-0.15cm}%
\includegraphics[width=.34\columnwidth]{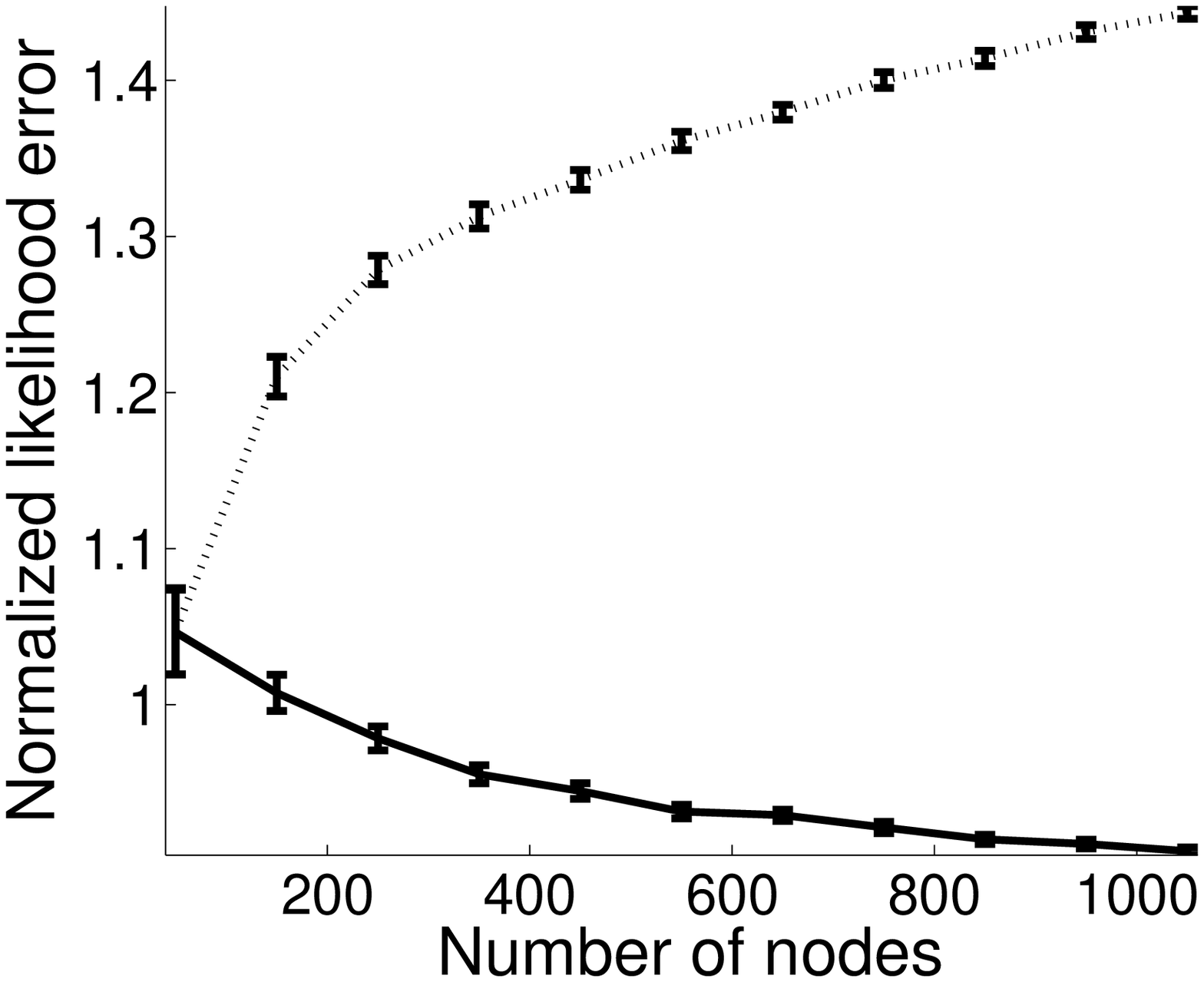}\hspace{-0.15cm}%
\includegraphics[width=.34\columnwidth]{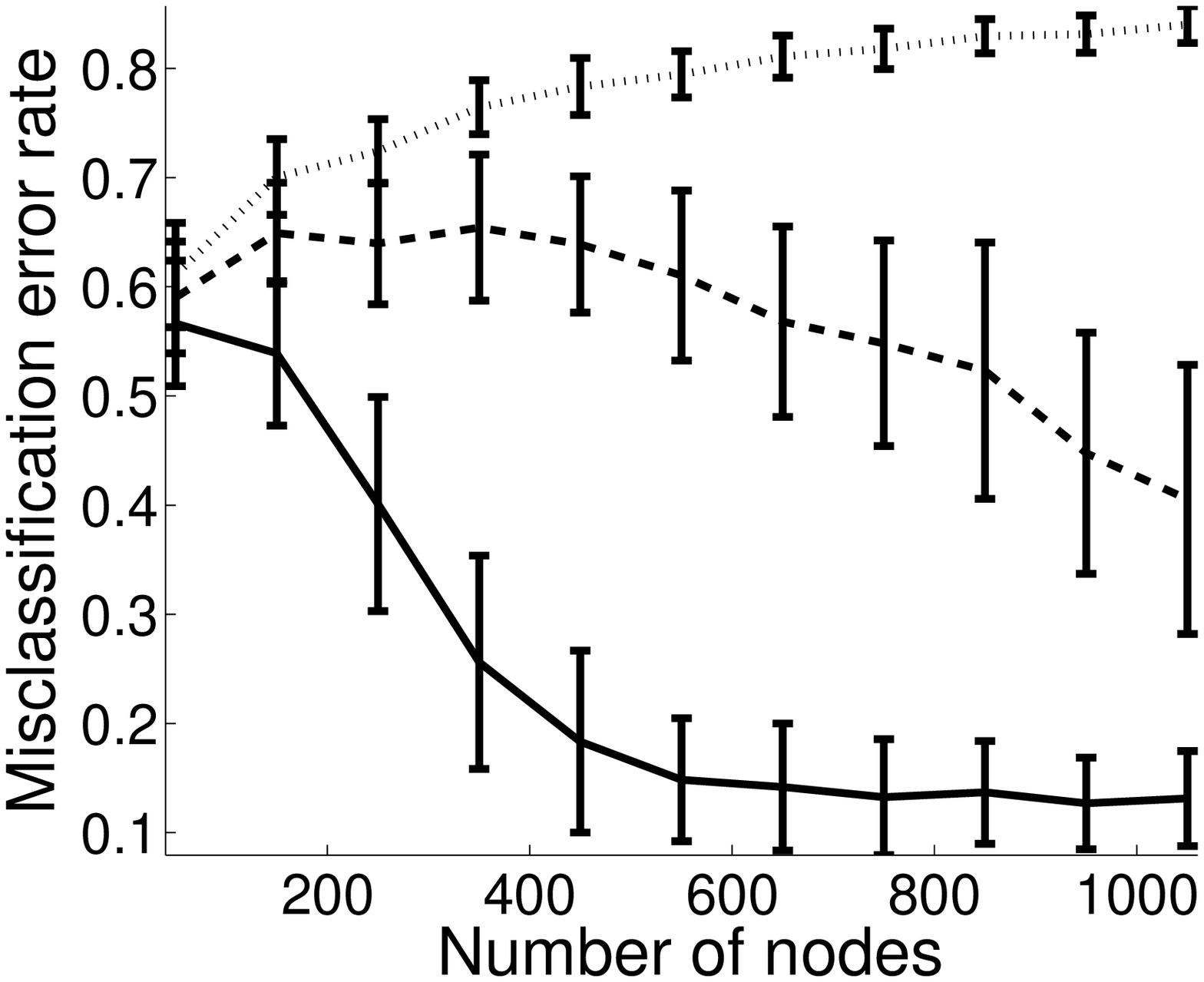}
\figurebox{-1.5pc}{35pc}{}[]
\caption{Simulation study results illustrating Theorems~\ref{thm:sumKLDivBnd}--\ref{thm:classConv}. Left: Likelihood error $|L(A;z) - \bar{L}_P(z)|/M$ as a function of network size $N$, shown for $M = N (\log N)^4$ with $K = N^{1/2}$.  Center: Same quantity for $M = N (\log N)^2$ with $K = N^{3/5}$ (dotted) and $K = N^{1/2}$ (solid).  Right: Error rate $N_{\mathrm{e}}(\hat{z})/N$ for $M = N (\log N)^2$ with $K = N^{1/2}$ and $\gamma = 4/5$ (dotted), $\gamma = 9/10$ (dashed), $\gamma = 1$ (solid)}
\label{fig:2}
\end{center}
\end{figure}

To test whether the assumptions of Theorem~\ref{thm:unifLikBnd} are necessary as well as sufficient to obtain convergence of $L(A;z)/M$ to $\bar{L}_P(z)/M$, blockmodels were next fitted to Erd\"{o}s-R\'{e}nyi networks of increasing size, for $N$ in the range 50--1050.  The corresponding normalized log-likelihood error $|L(A;z) - \bar{L}_P(z)|/M$ for different rates of growth in the expected number of edges $M$ and the number of fitted classes $K$ is shown in Fig.~\ref{fig:2}.  Observe from the leftmost panel that when $M= N (\log N)^4$ and $K = N^{1/2}$, as prescribed by the theorem, this error decreases in $N$. If the edge density is reduced to $M/N = (\log N)^2$, we observe in the center panel convergence when $K = N^{1/2}$ and divergence when $K = N^{3/5}$. This suggests that the error as a function of $K$ follows Theorem \ref{thm:unifLikBnd} closely, but that the network can be somewhat more sparse than it requires.

To test the conditions of Theorem~\ref{thm:classConv}, blockmodels with parameters $(\bar{z}, \bar{\theta})$ and increasing class size $K$ were used to generate data, and corresponding node misclassification error rates $N_{\mathrm{e}}(z)/N$ were recorded as a function of correctly specified $K$-class blockmodel fitting.  Model parameter $\bar{z}$ was chosen to yield equally-sized blocks, so as to meet identifiability condition (i) of Theorem~\ref{thm:classConv}.  Parameter $\bar{\theta} = \alpha I +\beta 1 1^{\mathrm{\scriptscriptstyle T}}$ was chosen to yield within-class and between-class success probabilities with the property that for any class pair $(a,b)$, the condition $D(\theta_{aa} \mid\mid (\theta_{aa}+\theta_{ab})/2) =  MK^\gamma/(20 N^2)$ was satisfied, with $\gamma\in \{4/5, 9/10, 1\}$; identifiability condition (ii) was thus met only in the $\gamma = 1$ case.  The rightmost panel of Fig.~\ref{fig:2} shows the fraction $N_{\mathrm{e}}(z)/N$ of misclassified nodes when $M = N (\log N)^2$ and $K = N^{1/2}$, corresponding to the setting in which convergence of $L(A;z)/M$ to $\bar{L}_P(z)/M$ was observed above; this fraction is seen to decay when $\gamma = 1$ or $9/10$, but to increase when $\gamma = 4/5$. This behavior conforms with Theorem~\ref{thm:classConv} and suggests that its identifiability conditions are close to being necessary as well as sufficient.

\section{Network data example}

\subsection{Facebook social network dataset}

To illustrate the use of our results in the fitting of $K$-class stochastic blockmodels to network data, we employed a publicly available social network dataset containing $N = 553$ undergraduate Facebook profiles from the California Institute of Technology (people.maths.ox.ac.uk/$\sim$porterm/data/facebook5.zip).  These profiles indicate whenever a pair of students have identified one another as friends, yielding a network of $11\,511$ edges and accompanying covariate information including gender, class year, and hall of residence.

\begin{figure}
\begin{center}
\hspace{-0.8cm}\includegraphics[width=0.38\columnwidth]{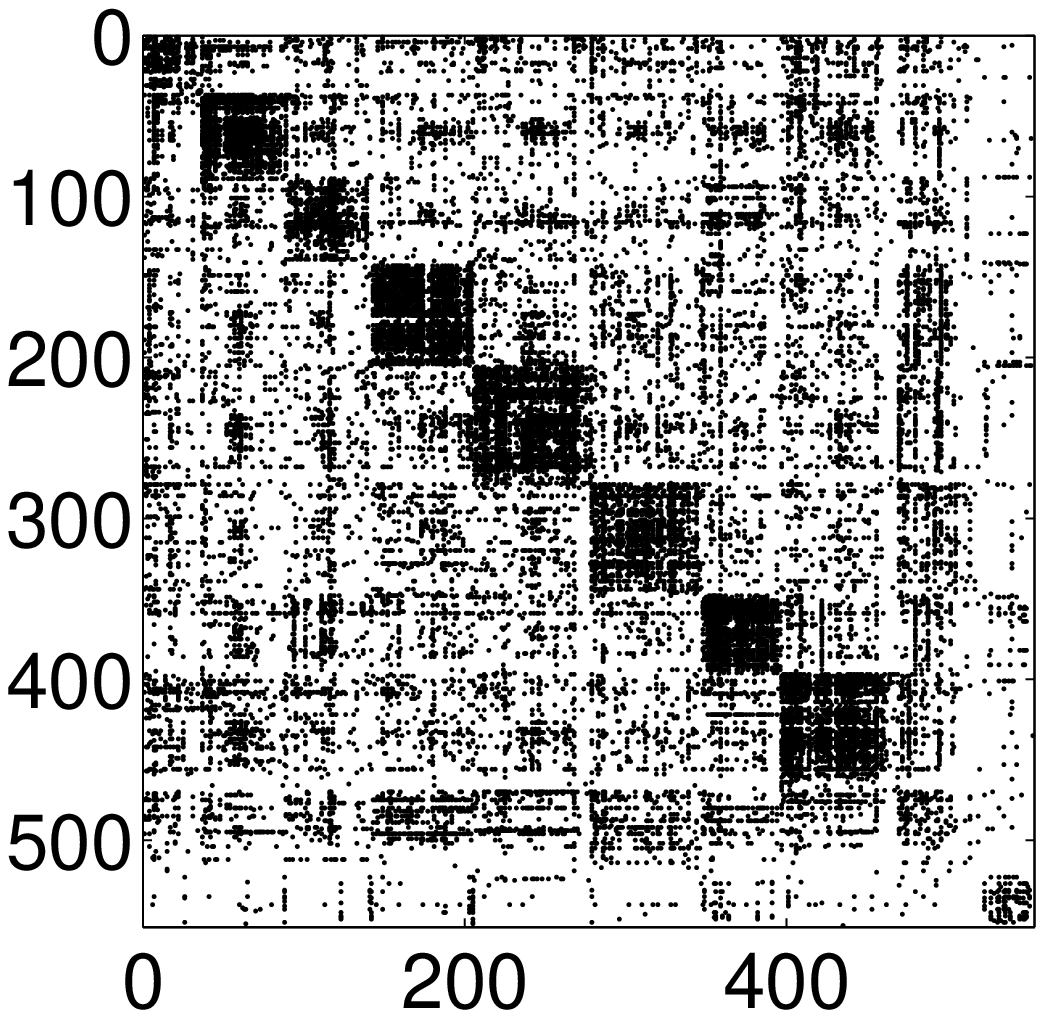}\hspace{-0.75cm}%
\includegraphics[width=0.37\columnwidth]{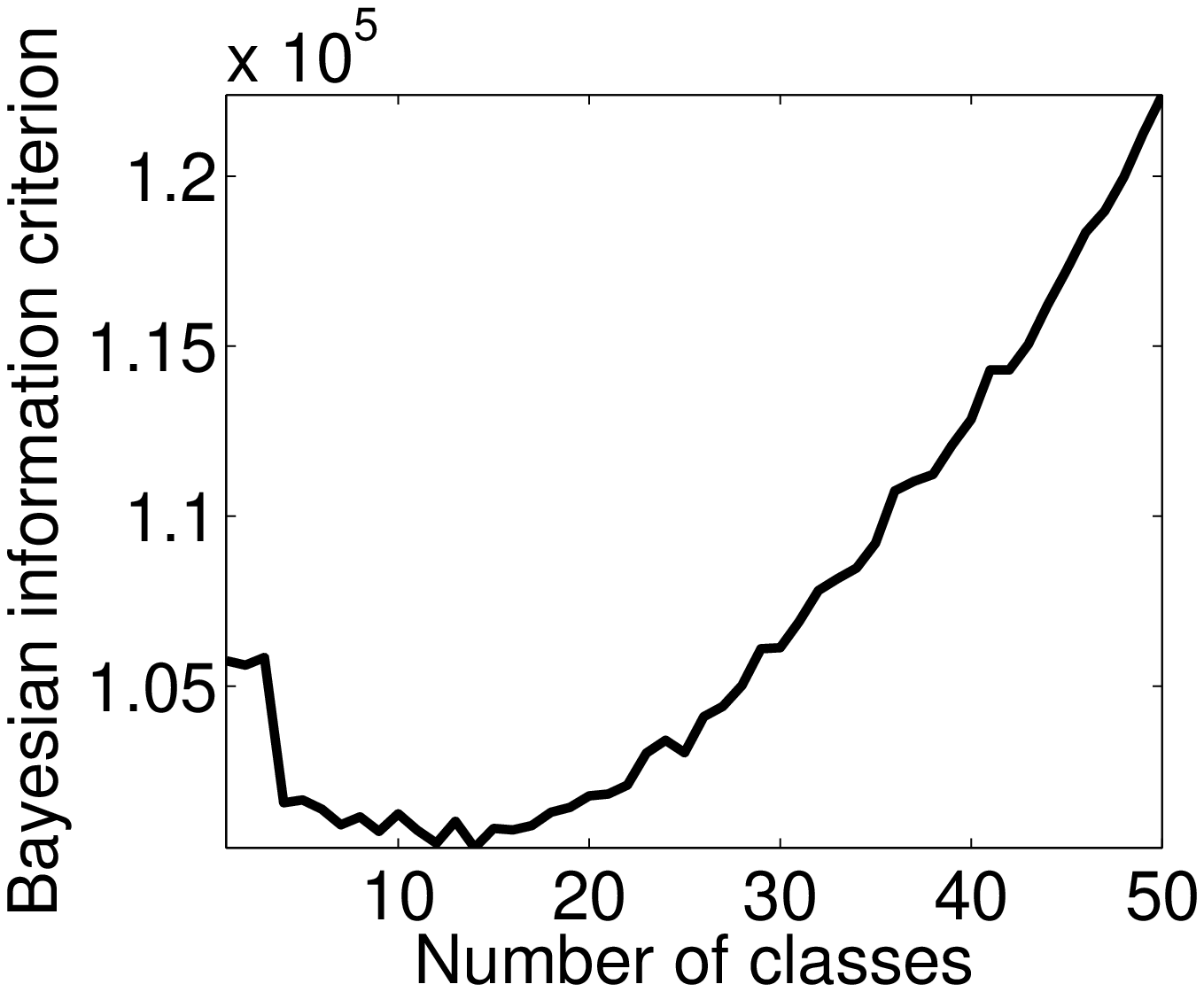}\hspace{-0.25cm}%
\includegraphics[width=0.37\columnwidth]{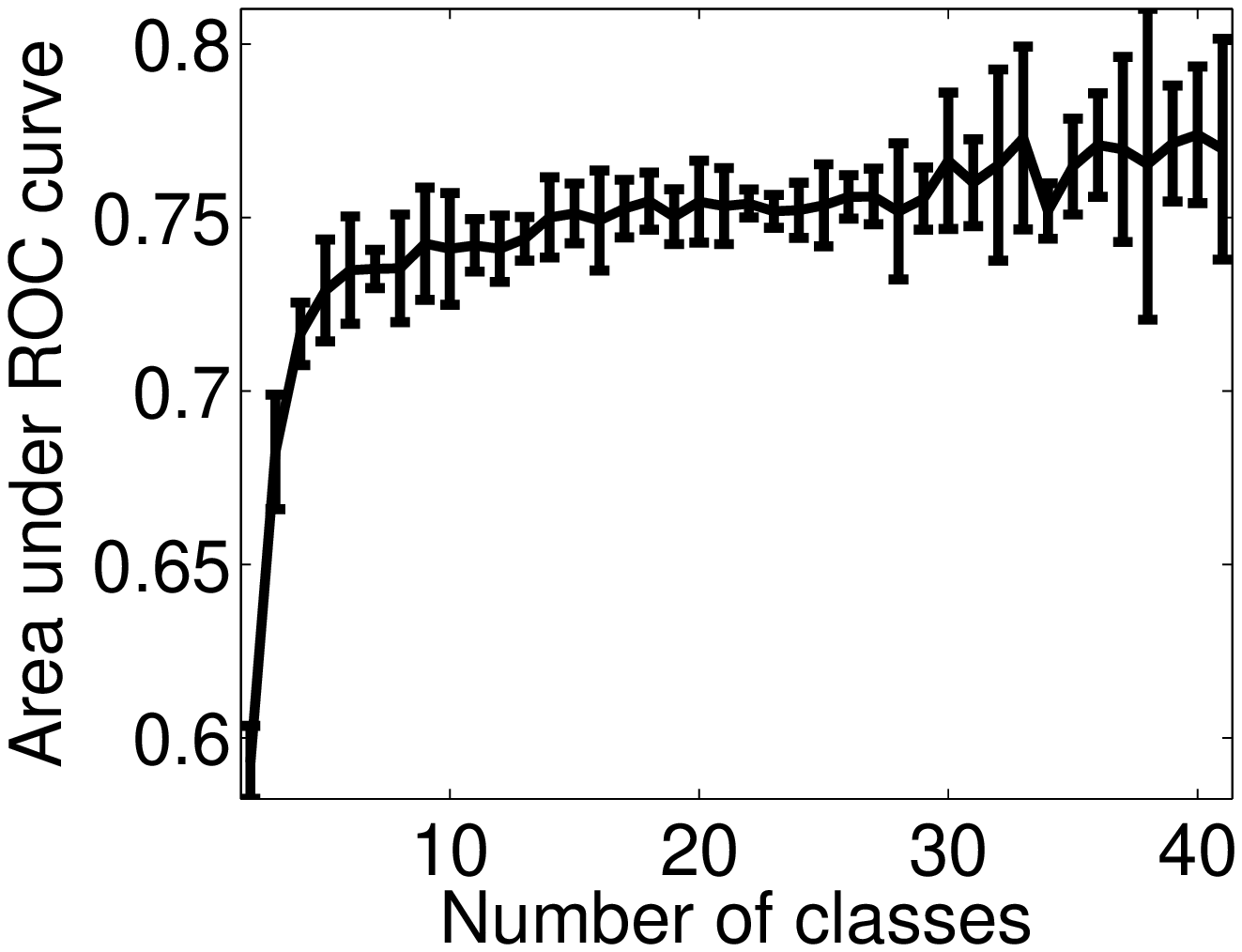}
\caption{Facebook social network dataset and its fitting statistics for varying number of blockmodel classes $K$.  Left: Adjacency data matrix of a network of Facebook undergraduate student profiles.  Center: Model order statistic for fitted logit blockmodels as a function of $K$.  Right: Out-of-sample prediction error as a function of $K$}
\label{fig:face1}
\end{center}
\end{figure}

\citet{traud2010community} applied community detection algorithms to this network, and compared their output to partitions based on categorical covariates such as those identified above.  They concludes that a grouping of students by residence hall was most similar to the best algorithmic grouping obtained, and thus that shared residence hall membership was the best predictor for the formation of community structure.  This structure is reflected in the leftmost panel of Fig.~\ref{fig:face1}, which shows the network adjacency structure under an ordering of students by residence hall.

\subsection{Logit blockmodel parameterization and fitting procedure}

Here we build on the results of \citet{traud2010community} by taking covariate information explicitly into account when fitting the Facebook dataset described above.  Specifically, by assuming only that links are independent Bernoulli variates and then employing confidence bounds to assess fitted blocks by way of parameter $\bar{\theta}^{(z)}$, we examine these data for residual community structure beyond that well explained by the covariates themselves.

Since the results of Theorems~\ref{thm:sumKLDivBnd} and~\ref{thm:unifLikBnd} hold uniformly over all choices of blockmodel membership vector $z$, we may select $z$ in any manner, including those that depend on covariates.  For this example, we determined an approximate maximum likelihood estimate $\hat{z}$ under a logit blockmodel that allows the direct incorporation of covariates.  The model is parameterized such that the log-odds ratio of an edge occurrence between nodes $i$ and $j$ is given by
\begin{equation}\label{eq:logitBlockmodel}
\log \frac{P_{ij}}{1-P_{ij}} = \tilde{\theta}_{z_i z_j} + x(i,j)^{\mathrm{\scriptscriptstyle T}} \beta \quad (i = 1, \ldots, N; j = i + 1, \ldots, N),
\end{equation}
where $x(i,j)$ a vector of covariates indicating shared group membership, and model parameters $(\tilde{\theta}, \beta, z)$ are estimated from the data.  Four categorical covariates were used: the three indicated above, plus an eight-category covariate indicating the range of the observed degree of each node; see \citet{karrer2011stochastic} for related discussion on this point.  Matrix $\tilde{\theta}$ is analogous to blockmodel parameter $\theta$, vector $z$ specifies the blockmodel class assignment, and vector $\beta$ was implemented here with sum-to-zero identifiability constraints.

Because exact maximization of the log-likelihood function $L(A;\tilde{\theta}, \beta, z)$ corresponding to~\eqref{eq:logitBlockmodel} is computationally intractable, we instead employed an approach that alternated between Markov chain Monte Carlo exploration of $z$ while holding $(\tilde{\theta}, \beta)$ constant, and optimization of $\tilde{\theta}$ and $\beta$ while holding $z$ constant.  We tested different initialization methods and observed that highest likelihoods were consistently produced by first fitting class assignment vector $z$.  This fitting procedure provides a means of estimating averages $\bar{\theta}^{(z)}$ over subsets of the set $\{P_{ij}\}_{i<j}$, under the assumption that the network data comprise independent $\operatorname{Bernoulli}(P_{ij})$ trials.
\begin{figure}
\begin{center}
\hspace{-0.8cm}%
\includegraphics[width=0.3\columnwidth]{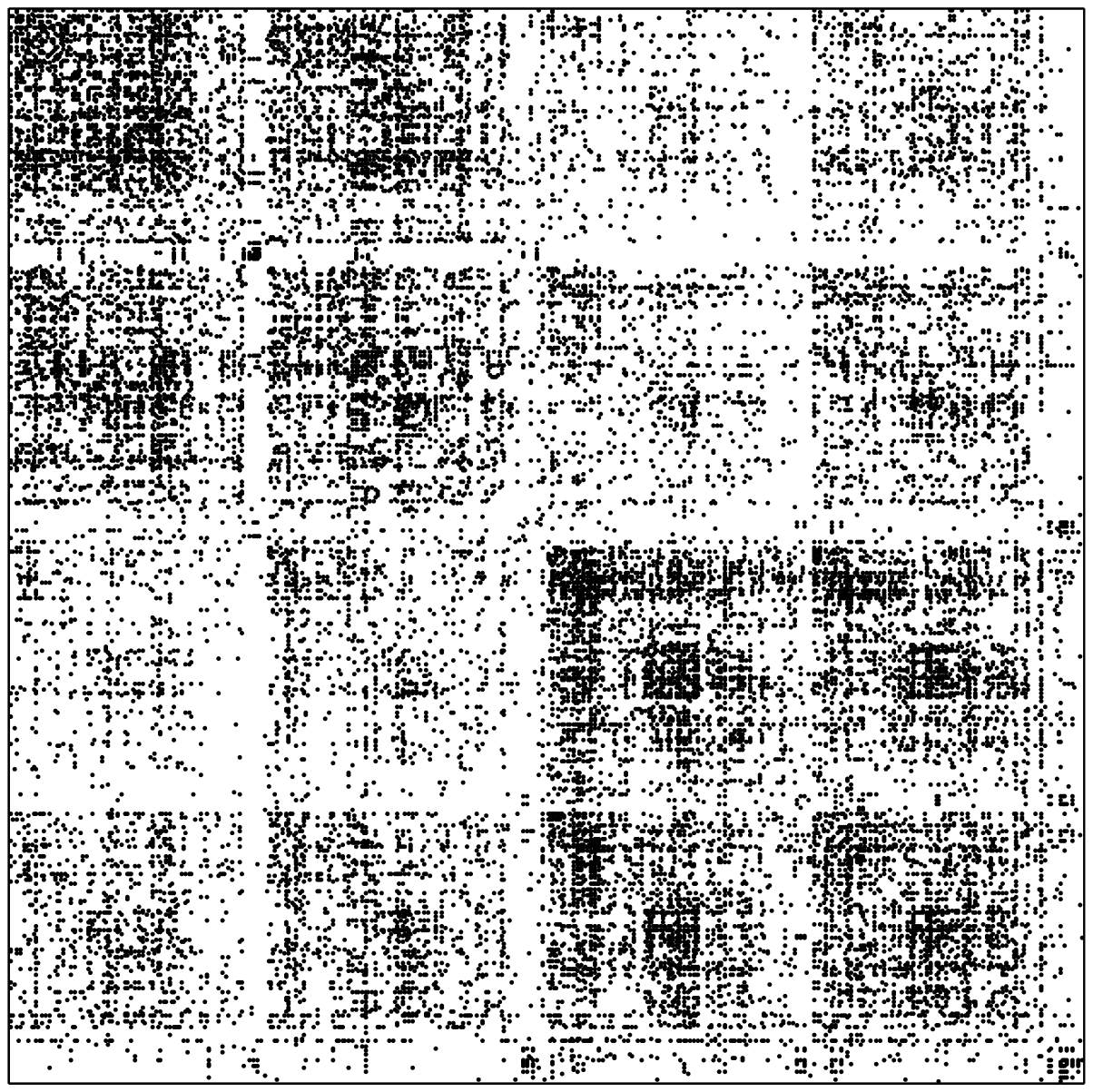}\hspace{-0.8cm}%
\includegraphics[width=0.3\columnwidth]{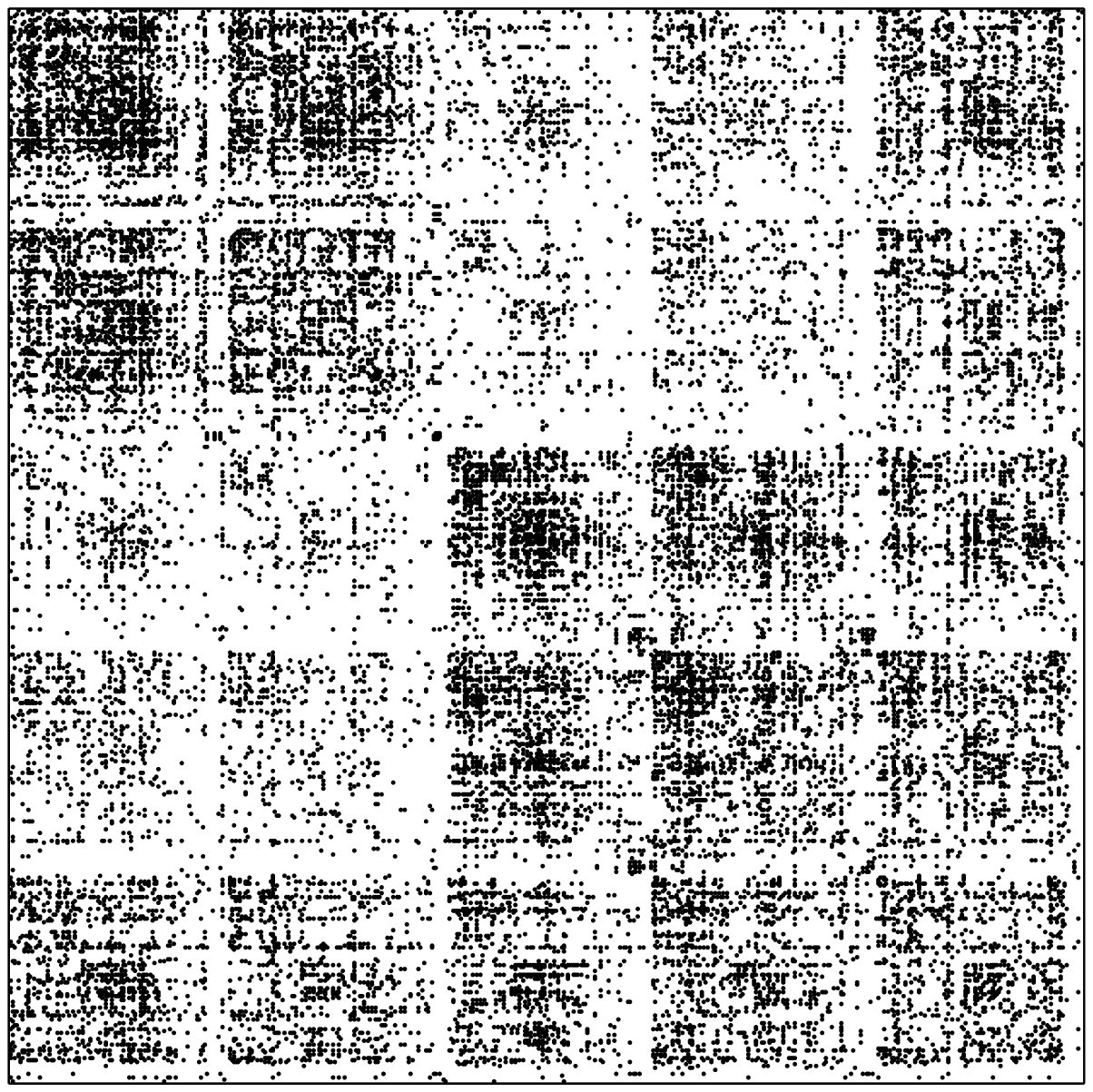}\hspace{-0.8cm}%
\includegraphics[width=0.3\columnwidth]{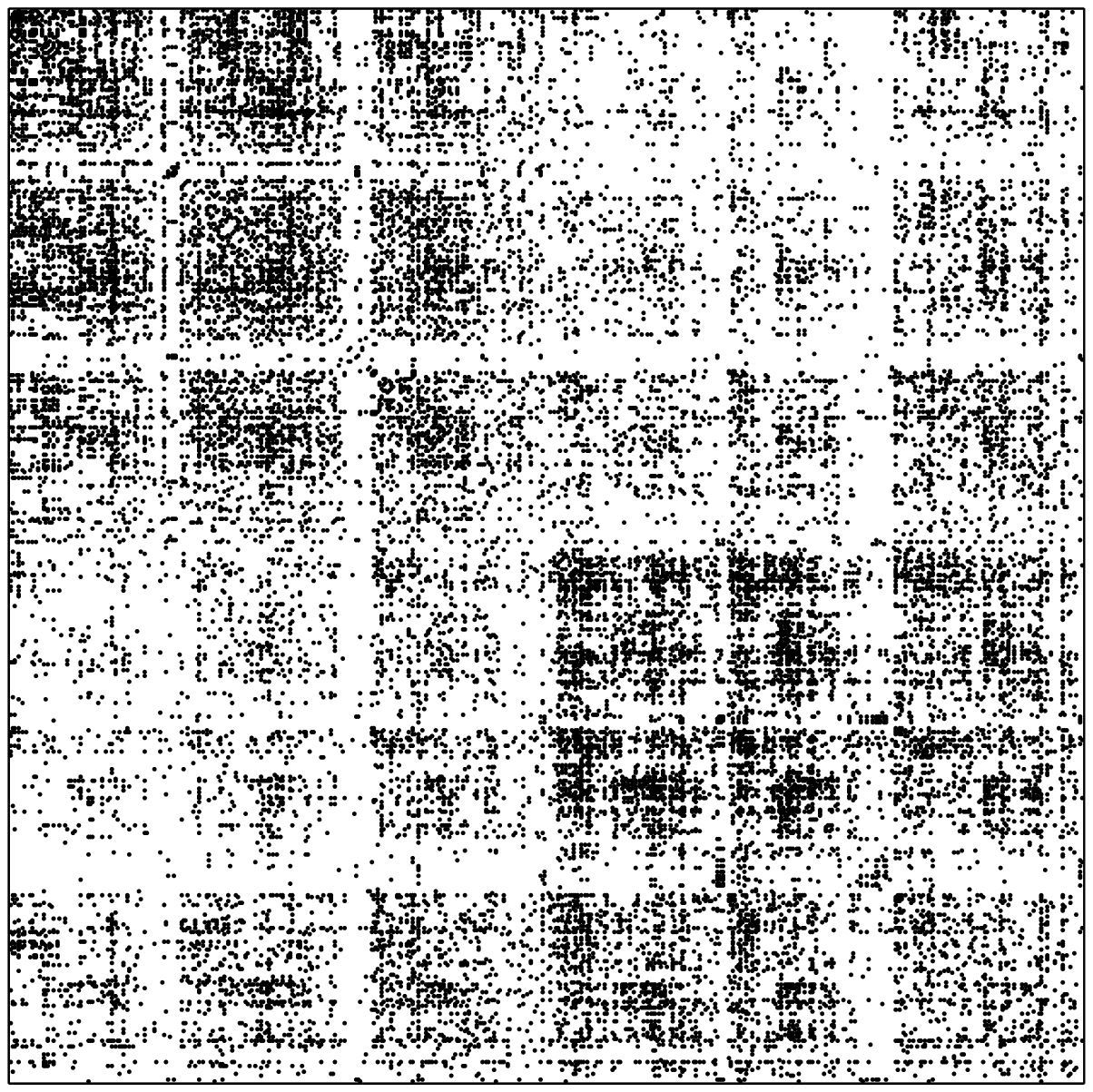}\hspace{-0.8cm}%
\includegraphics[width=0.3\columnwidth]{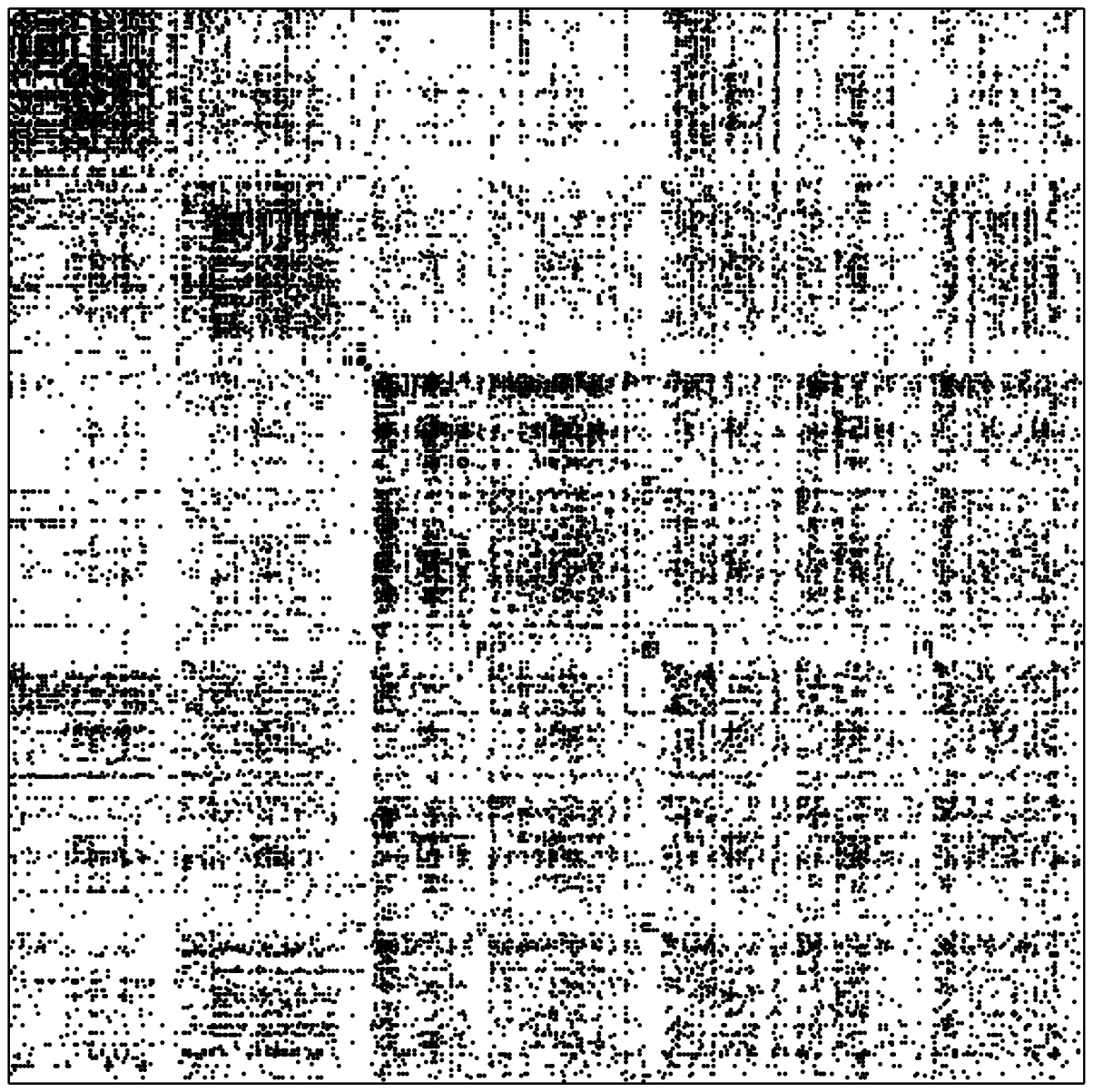}\hspace{-0.8cm}%
\hspace{-0.cm}%
\includegraphics[width=0.25\columnwidth]{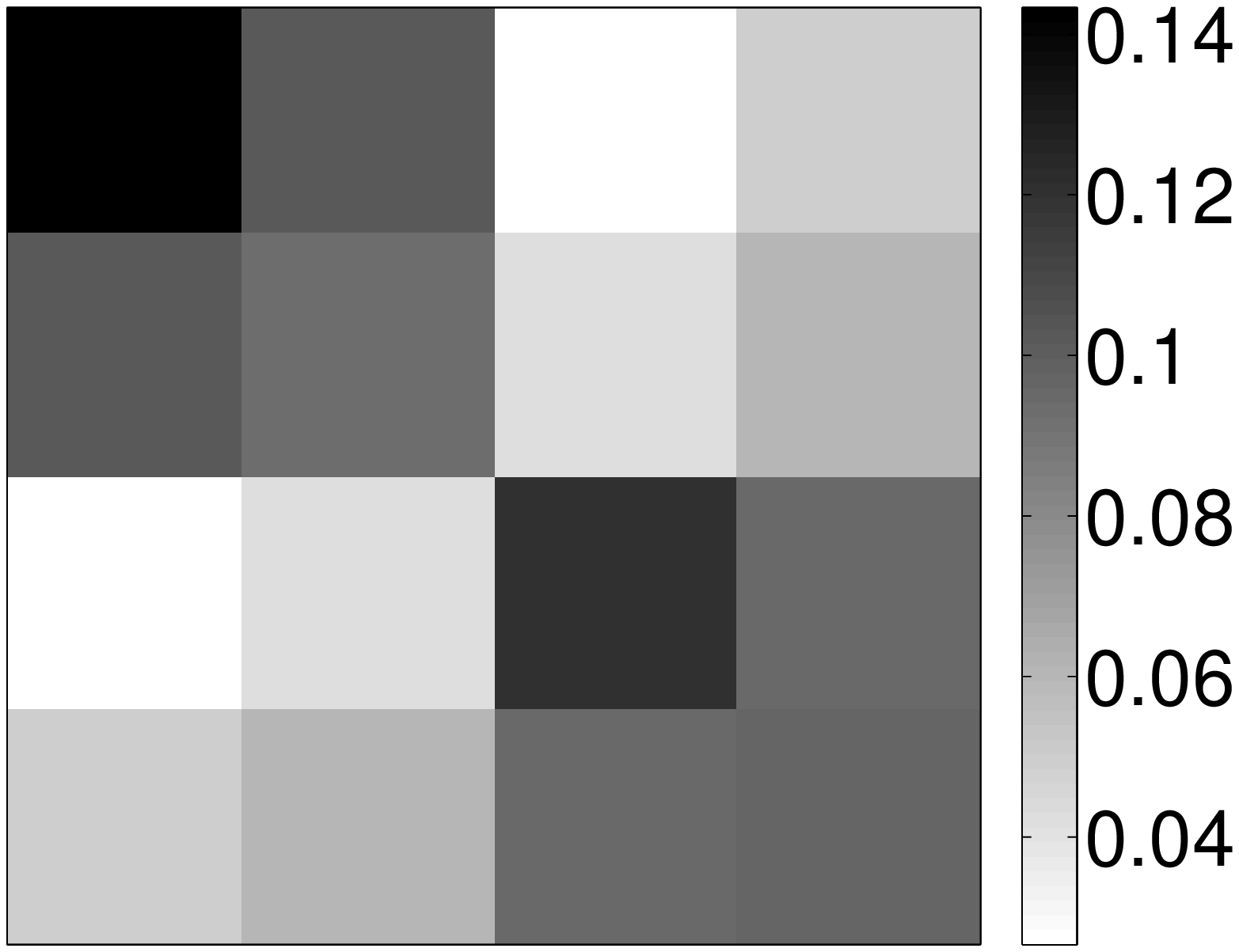}\hspace{-0.05cm}%
\includegraphics[width=0.25\columnwidth]{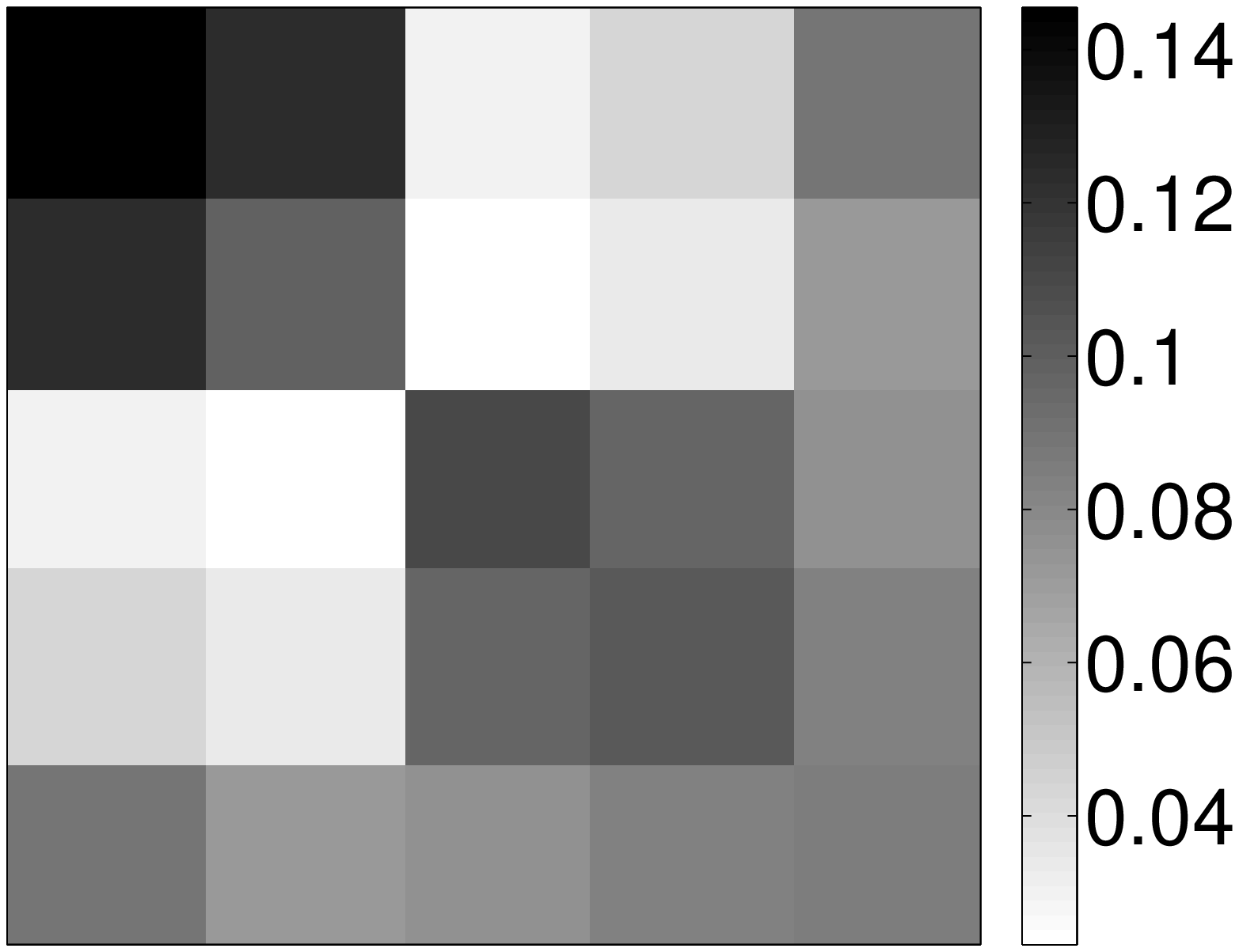}\hspace{-0.05cm}%
\includegraphics[width=0.25\columnwidth]{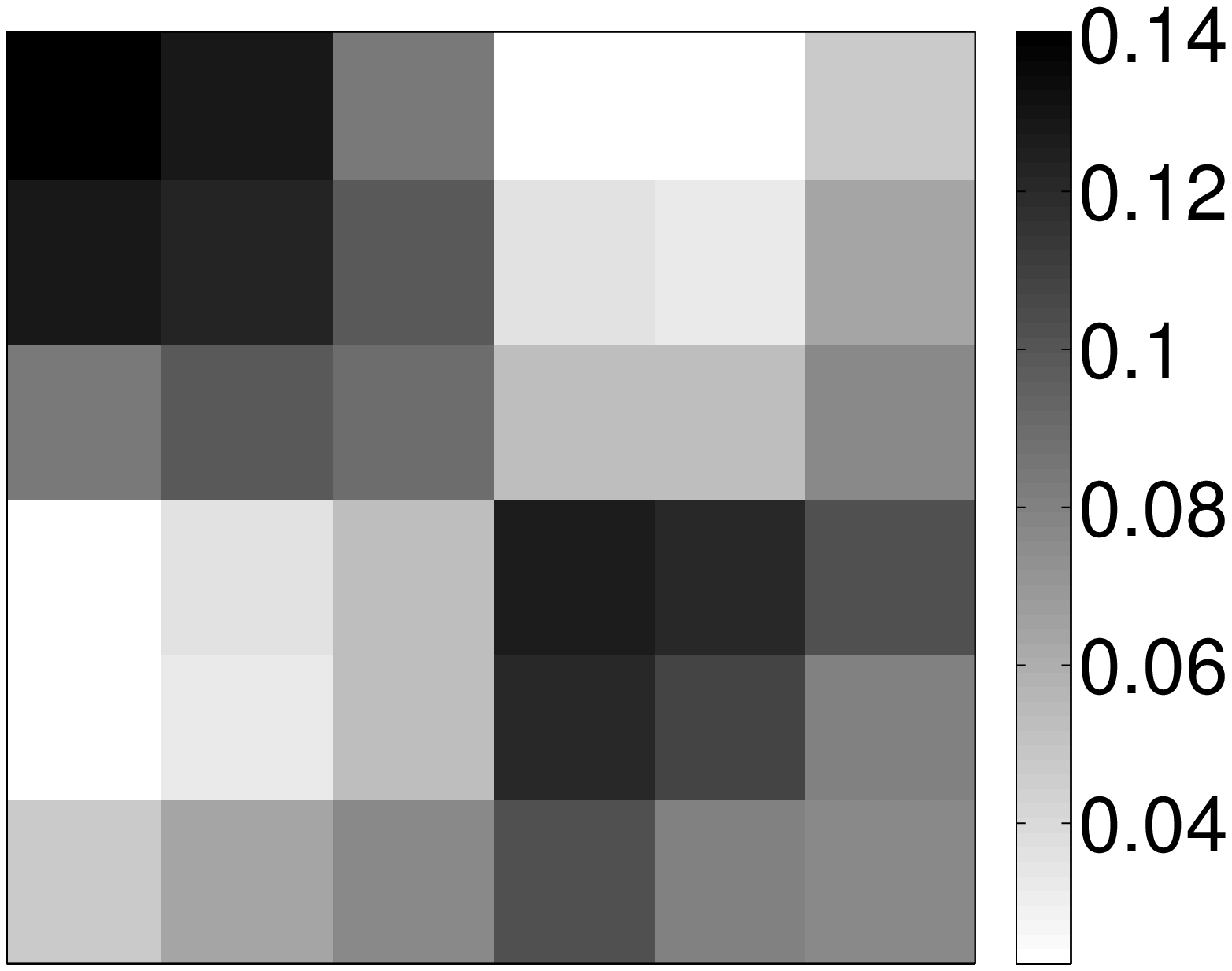}\hspace{-0.05cm}%
\includegraphics[width=0.25\columnwidth]{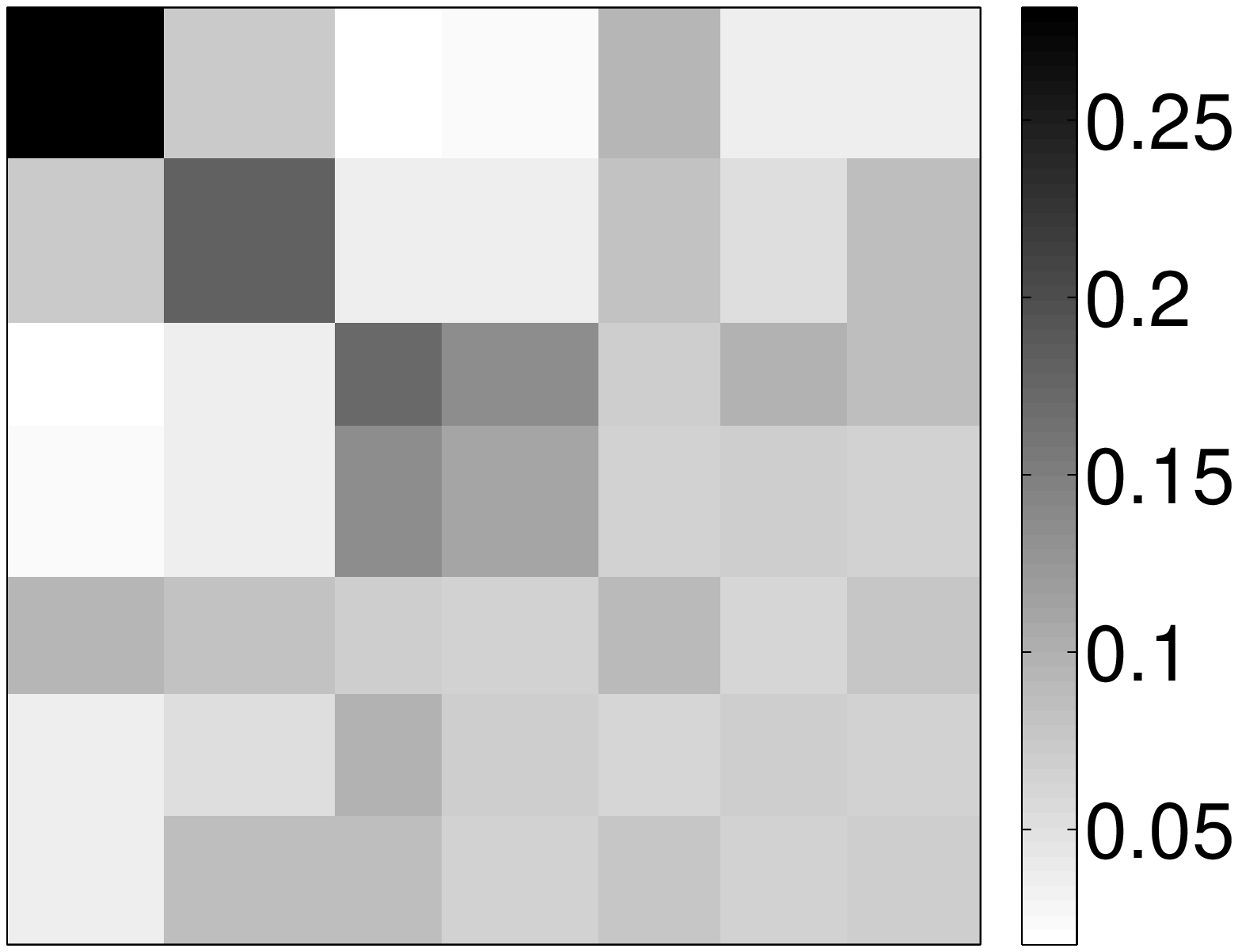}\hspace{-0.05cm}%
\includegraphics[width=0.33\columnwidth]{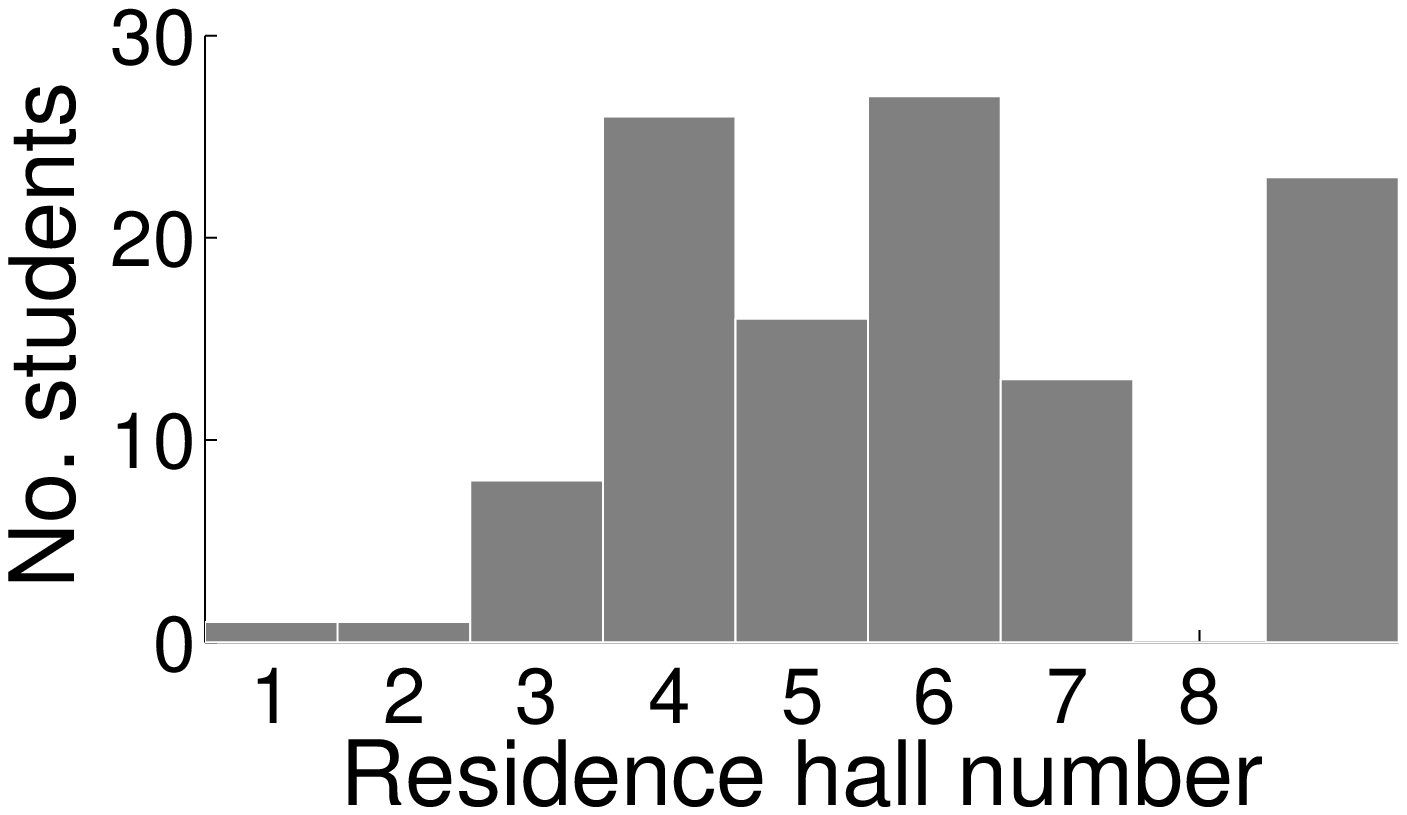}
\includegraphics[width=0.33\columnwidth]{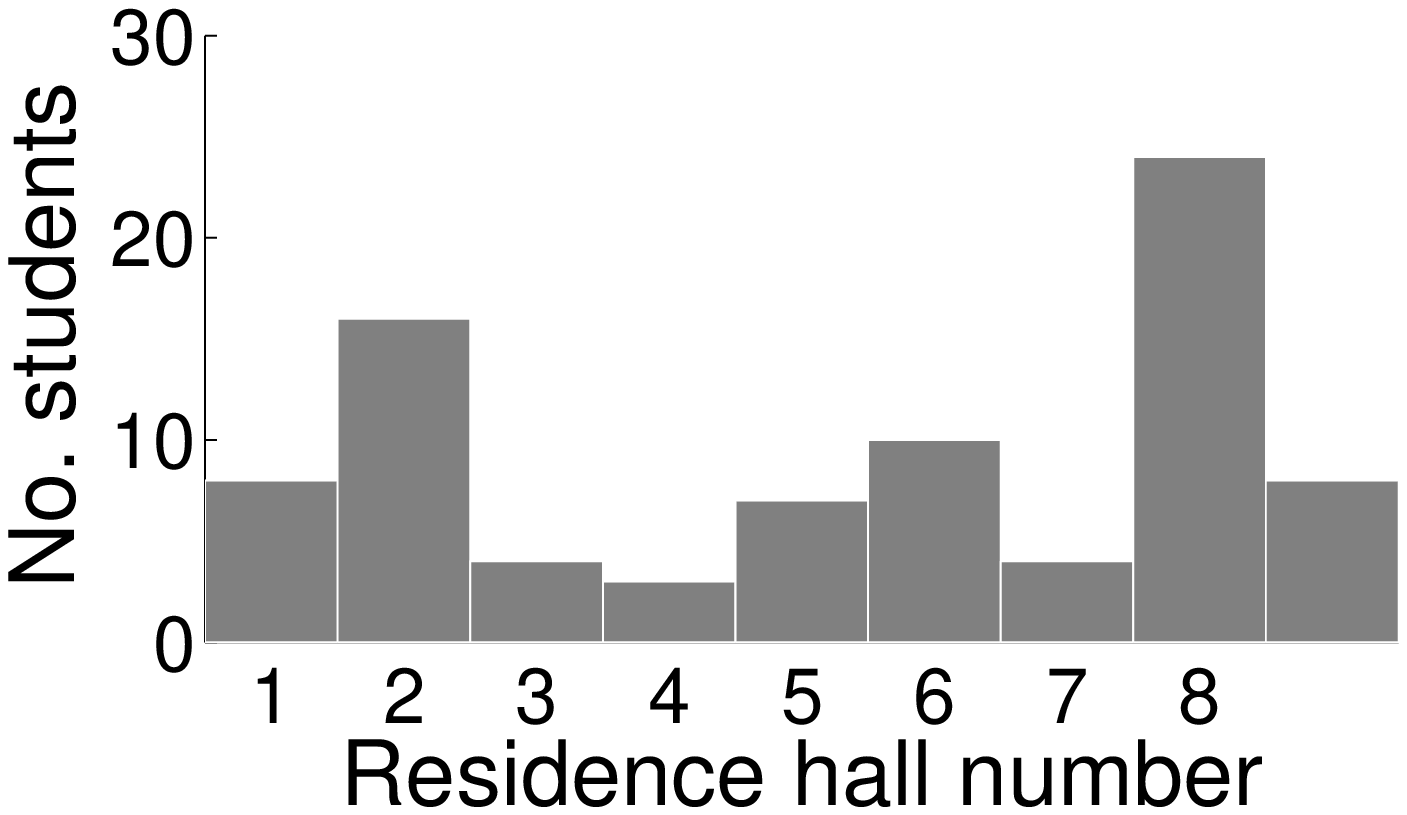}
\caption{Results of logit blockmodel fitting to the data of Fig.~\ref{fig:face1} for each of $K \in \{4,5,6,7\}$ classes.  Top row: Adjacency structure of the data, permuted to show block assignments for $K \in \{4,5,6,7\}$.  Second row: Corresponding estimates $\hat{\theta}$, with Kullback--Leibler divergence bounds 0$\cdot$0057, 0$\cdot$0067, 0$\cdot$0077, and 0$\cdot$0086.  Bottom row: Residence hall assignments of students whose grouping remained constant over these four values of $K$}
\label{fig:face8}
\end{center}
\end{figure}

\subsection{Data analysis}

We fitted the logit blockmodel of~\eqref{eq:logitBlockmodel} for values of $K$ ranging from $1$ to $50$ using the stochastic maximization procedure described in the preceding paragraph, and gauged model order by the Bayesian information criterion and out-of-sample prediction using five-fold cross validation, shown respectively in the center and rightmost panels of Fig.~\ref{fig:face1}.  These plots suggest a relatively low model order, beginning around $K=4$.  The corresponding 95\% confidence bounds on the divergence of $\hat{\theta}^{(z)}$ from $\bar{\theta}^{(z)}$ provided by Theorem~\ref{thm:sumKLDivBnd} also yield small values for $K$ in the range 4--7: for example, when $K=5$, the normalized sum of Kullback--Leibler divergences ${N \choose 2} {}^{-1} \sum_{a\leq b} n_{ab} D(\hat{\theta}_{ab} \mid\mid \bar{\theta}_{ab})$ is bounded by 0$\cdot$0067.  Corresponding normalized root-mean-square error bounds over this range of $K$ are approximately one order of magnitude larger.

We then examined approximate maximum likelihood estimates of $z$ for $K$ in the range 4--7, as shown in the top two rows of Fig.~\ref{fig:face8}; larger values of $K$ also reveal block structure, but exhibit correspondingly larger confidence bound evaluations.  The permuted adjacency structures under each estimated class assignment $\hat{z}$ are shown in the top row, along with the corresponding values of $\hat{\theta}$ below in the second row.  The structure of $\hat{\theta}$ over this range of $K$ suggests that after covariates are taken into account, it is possible to identify a subset of students who divide naturally into two residual ``meta-groups'' that interact less frequently with one another in comparison to the remaining subjects in the dataset; the precision of the corresponding estimates $\hat{\theta}$ can be quantified by Theorem~\ref{thm:sumKLDivBnd}, as in the caption of Fig.~\ref{fig:face8}.

As $K$ increases, these groups become more tightly concentrated, as extra blocks absorb students whose connections are more evenly distributed.  While the exact membership of each group varied over $K$, in part due to stochasticity in the fitting algorithm employed, we observed 199 students whose meta-group membership remained constant.  The bottom row of Fig.~\ref{fig:face8} shows the 8 residence halls identified for these sets of students, with the ninth category indicating unreported; observe that the effect of residence hall is still visible in that the left-hand grouping has more students in halls 4--7, while the right-hand grouping has more students in halls 1, 2, and 8.

\section*{Acknowledgement}

Work supported in part by the National Science Foundation, National Institute of Health, Army Research Office and the Office of Naval Research, U.S.A. Additional funding provided by the Harvard Medical School's Milton Fund.

\appendix
\section*{Appendix}

\subsection*{Proofs of Theorems~\ref{thm:sumKLDivBnd} and~\ref{thm:unifLikBnd}}

\begin{proof}[of Theorem~\ref{thm:sumKLDivBnd}]
To begin, observe that for any fixed class assignment $z$, every $\hat{\theta}_{ab}$ is a sum of $n_{ab}$ independent Bernoulli random variables, with corresponding mean $\bar{\theta}_{ab}$.  A Chernoff bound~\citep{dubhashi2009concentration} shows
\begin{align*}
\operatorname{pr}(\hat{\theta}_{ab} & \geq \bar{\theta}_{ab}+t ) \leq e^{-n_{ab} D(\bar{\theta}_{ab}+t \mid\mid \bar{\theta}_{ab})}, \quad 0 < t \leq 1 - \bar{\theta}_{ab}\\
\operatorname{pr}(\hat{\theta}_{ab} & \leq \bar{\theta}_{ab}-t ) \leq e^{-n_{ab} D(\bar{\theta}_{ab}-t \mid\mid \bar{\theta}_{ab})},
\quad 0 < t \leq \bar{\theta}_{ab}.
\end{align*}

Since these bounds also hold respectively for $\operatorname{pr}(\hat{\theta}_{ab} = \bar{\theta}_{ab} \pm t )$, we may bound the probability of any given realization $\vartheta \in \{0,1/n_{ab},\ldots,1\}$ of $\hat{\theta}_{ab}$ in terms of the Kullback--Leibler divergence of $\bar{\theta}_{ab}$ from $\vartheta$:
$$
\operatorname{pr}(\hat{\theta}_{ab} = \vartheta) \leq e^{-n_{ab} D(\vartheta \mid\mid \bar{\theta}_{ab})}.
$$
By independence of the $\{A_{ij}\}_{i< j}$, this implies a corresponding bound on the probability of any $\hat{\theta}$:
\begin{equation}\label{eq:thetaHatBnd}
\operatorname{pr}(\hat{\theta}) \leq \exp\left\{-\textstyle \sum_{a \leq b} n_{ab} D(\hat{\theta}_{ab} \mid\mid \bar{\theta}_{ab})\right\}.
\end{equation}

Now, let $\widehat{\Theta}$ denote the range of $\hat{\theta}$ for fixed $z$, and observe that since each of the $\binom{K+1}{2}$ lower-diagonal entries $\{\hat{\theta}_{ab}\}_{a \leq b}$ of $\hat{\theta}$ can independently take on $n_{ab}+1$ distinct values, we have that $|\widehat{\Theta}| = \prod_{a\leq b} (n_{ab}+1)$.  Subject to the constraint that $\sum_{a\leq b} n_{ab} = \binom{N}{2}$, we see that this quantity is maximized when $n_{ab} = \binom{N}{2} / \binom{K+1}{2}$ for all $a \leq b$, and hence
\begin{equation}\label{eq:thetaCardBnd}
|\widehat{\Theta}| \leq \left[\textstyle \binom{N}{2}/\binom{K+1}{2}+1\right]^{\binom{K+1}{2}} \! < \left(N^2/K^2 + 1\right)^{\frac{K^2+K}{2}} \! < \left(N/K + 1\right)^{K^2+K} \! .
\end{equation}

Now consider the event that $\sum_{a \leq b} n_{ab} D(\hat{\theta}_{ab} \mid\mid \bar{\theta}_{ab})$ is at least as large as some $\epsilon > 0$; the probability of this event is given by $\operatorname{pr}(\widehat{\Theta}_{\epsilon})$ for
\begin{equation}\label{eq:thetaEpsDef}
\widehat{\Theta}_{\epsilon} = \left\{ \hat{\theta} \in \widehat{\Theta} : \textstyle \sum_{a \leq b} n_{ab} D(\hat{\theta}_{ab} \mid\mid \bar{\theta}_{ab}) \geq \epsilon \right\}.
\end{equation}
Since $\sum_{a \leq b} n_{ab} D(\hat{\theta}_{ab} \mid\mid \bar{\theta}_{ab}) \geq \epsilon$ for all $\hat{\theta} \in \widehat{\Theta}_{\epsilon}$, we have from~\eqref{eq:thetaHatBnd} and~\eqref{eq:thetaEpsDef} that
$$
\operatorname{pr}(\widehat{\Theta}_{\epsilon}) = \sum_{\hat{\theta} \in \widehat{\Theta}_{\epsilon}} \operatorname{pr}(\hat{\theta}) \leq \sum_{\hat{\theta} \in \widehat{\Theta}_{\epsilon}} e^{-\sum_{a \leq b} n_{ab} D(\hat{\theta}_{ab} \mid\mid \bar{\theta}_{ab})}
\leq \sum_{\hat{\theta} \in \widehat{\Theta}_{\epsilon}} e^{-\epsilon} = |\widehat{\Theta}_{\epsilon}| e^{-\epsilon},
$$
and since $|\widehat{\Theta}_{\epsilon}| \leq |\widehat{\Theta}|$, we may use~\eqref{eq:thetaCardBnd} to obtain, for fixed class assignment $z$,
\begin{equation}\label{eq:divBndFixedz}
\operatorname{pr}\left\{\textstyle \sum_{a\leq b} n_{ab} D(\hat{\theta} \mid\mid \bar{\theta}) \geq \epsilon \right\} < \left( N/K+1\right)^{K^2+K} e^{-\epsilon}.
\end{equation}
Appealing to a union bound over all $K^N$ possible class assignments and setting $\epsilon = \log [K^N \left(N/K+1\right)^{K^2+K}/\delta]$ then yields the claimed result.
\end{proof}

\begin{proof}[of Theorem~\ref{thm:unifLikBnd}]
By Lemma~\ref{lem:LbarLDiff}, the difference $L(A;z) - \bar{L}_P(z)$ can be expressed for any fixed class assignment $z$ as $\sum_{a\leq b} n_{ab} D(\hat{\theta}_{ab} \mid\mid \bar{\theta}_{ab}) + X - E(X)$,
where the first term satisfies the deviation bound of~\eqref{eq:divBndFixedz}, and $X = \sum_{i<j} A_{ij} \log \{\bar{\theta}_{z_i z_j} / (1-\bar{\theta}_{z_iz_j})\}$ comprises a weighted sum of independent $\operatorname{Bernoulli}(P_{ij})$ random variables.

To bound the quantity $|X - E(X)|$, observe that since by assumption $N^{-2} \leq P_{ij} \leq 1-N^{-2}$, the same is true for each corresponding average $\bar{\theta}_{z_i z_j}$.  As a result, the random variables $X_{ij} = A_{ij} \log \{\bar{\theta}_{z_i z_j} / (1-\bar{\theta}_{z_iz_j})\}$ comprising $X$ are each bounded in magnitude by $C = 2 \log N$. This allows us to apply a Bernstein inequality for sums of bounded independent random variables due to~\citet[Theorems~2.8 and 2.9, p. 27]{chung2006complex}, which states that for any $\epsilon > 0$,
\begin{equation}
\label{eq:Xdev}
\operatorname{pr}\{|X - E(X) | \geq \epsilon \}  \leq 2 \exp\left\{-\frac{\epsilon^2}{2\sum_{i<j} E(X_{ij}^2) + (2/3)\epsilon C}\right\}.
\end{equation}

Finally, observe that since the event $|L(A;z) - \bar{L}_P(z)| > 2\epsilon M$ implies either the event $\sum_{a\leq b} n_{ab} D(\hat{\theta}_{ab} \mid\mid \bar{\theta}_{ab}) \geq \epsilon M$ or the event $|X - E(X) | \geq \epsilon M$, we have for fixed assignment $z$ that
\begin{align*}
\operatorname{pr}\{|L(A;z) - \bar{L}_P(z) & \geq 2\epsilon M\} \leq \textstyle \operatorname{pr}\Big[\Big\{\sum_{a\leq b} n_{ab} D(\hat{\theta}_{ab} \mid\mid \bar{\theta}_{ab}) \geq \epsilon M \Big\} \cup \Big\{|X - E(X) | \geq \epsilon M\Big\}\Big].
\end{align*}
Summing the right-hand sides of~\eqref{eq:divBndFixedz} and~\eqref{eq:Xdev}, and then over all $K^N$ possible assignments, yields
\begin{multline*}
\operatorname{pr}\{\max_z |L(A;z) - \bar{L}_P(z)| \geq 2\epsilon M\} \leq
\exp\left\{K\log N + (K^2+K)\log (N/K + 1) - \epsilon M\right\} \\ + 2\exp\left\{K\log N -\frac{\epsilon^2 M}{8\log^2 N+ (4/3)\epsilon\log N}\right\},
\end{multline*}
where we have used the fact that $\sum_{i<j} E(X_{ij}^2) \leq 4M \log^2 N$ in~\eqref{eq:Xdev}.  It follows directly that if $K = \mathcal{O}(N^{1/2})$ and $M = \omega(N (\log N)^{3+\delta})$, then $\lim_{N\rightarrow\infty} \operatorname{pr}\{\max_z |L(A;z) - \bar{L}_P(z)| / M \geq \epsilon\} = 0$ for every fixed $\epsilon > 0$ as claimed.
\end{proof}

\subsection*{Proof of Theorem~\ref{thm:classConv}}

\begin{proof}[of Theorem~\ref{thm:classConv}]
To begin, note that Theorem \ref{thm:unifLikBnd} holds uniformly in $z$, and thus implies that
\begin{equation*}
|\bar{L}_P(\bar{z}) - L(A;\bar{z})| + |\bar{L}_P(\hat{z}) - L(A;\hat{z})| = o_P(M).
\end{equation*}
Since $\hat{z}$ is the maximum-likelihood estimate of class assignment $\bar{z}$, we know that $L(A;\hat{z}) \geq L(A;\bar{z})$, implying that $L(A;\hat{z}) = L(A;\bar{z}) + \delta$ for some $\delta \geq 0$.  Thus, by the triangle inequality,
\begin{equation*}
|\bar{L}_P(\bar{z}) - \bar{L}_P(\hat{z}) + \delta| \leq |\bar{L}_P(\bar{z}) - L(A;\bar{z})| + |\bar{L}_P(\hat{z}) - (L(A;\bar{z}) + \delta)| = o_P(M),
\end{equation*}
and since $\bar{L}_P(\bar{z}) \geq \bar{L}_P(\hat{z})$ under any blockmodel with parameter $\bar{z}$, we have $\bar{L}_P(\bar{z}) - \bar{L}_P(\hat{z}) = o_P(M)$.

Under conditions (i) and (ii) of Theorem~\ref{thm:classConv}, we will now show that also
\begin{equation} \label{eq:NerrBnd}
\bar{L}_P(\bar{z}) - \bar{L}_P(\hat{z}) = \frac{N_{\mathrm{e}}(\hat{z})}{N} \, \Omega (M),
\end{equation}
holds for \emph{every} realization of $\hat{z}$,
thus implying that $N_{\mathrm{e}}(\hat{z}) = o_P(N)$ and proving the theorem.

To show~\eqref{eq:NerrBnd}, first observe that any blockmodel class assignment vector $z$ induces a corresponding partition of the set $\{P_{ij}\}_{i<j}$ according to $(i,j) \mapsto (z_i,z_j)$.  Formally, $z$ partitions $\{P_{ij}\}_{i<j}$ into $L$ subsets $(S_1, \ldots, S_L)$ via the mapping
$$
\zeta_{ij}: (i=1,\ldots,N; j=i+1,\ldots,N) \to (l=1,\ldots,L).
$$
This partition is separable in the sense that there exists a bijection between $\{1,\ldots,L\}$ and the upper triangular portion of blockmodel parameter $\theta$, such that we write $\theta_{\zeta_{ij}} = \theta_{z_i z_j}$ for membership vector $z$.  More generally, for \emph{any} partition $\Pi$ of $\{P_{ij}\}_{i<j}$, we may define $\bar{\theta}_l = |S_l|^{-1} \sum_{i < j} P_{ij} \, 1\{ P_{ij} \in S_l \}$ as the arithmetic average over all $P_{ij}$ in the subset $S_l$ indexed by $\zeta_{ij} = l$.  Thus we may also define
$$
\bar{L}_P^*(\Pi) = \sum_{i<j} \left\{ P_{ij} \log \bar{\theta}_{\zeta_{ij}} + (1-P_{ij}) \log (1-\bar{\theta}_{\zeta_{ij}}) \right\},
$$
so that $\bar{L}_P^*$ and $\bar{L}_P$ coincide on partitions corresponding to admissible blockmodel assignments $z$.

The establishment of~\eqref{eq:NerrBnd} proceeds in three steps: first, we construct and analyze a refinement of the partition $\Pi^z$ induced by any blockmodel assignment vector $z$ in terms of its error $N_{\mathrm{e}}(z)$; then, we show that refinements increase $\bar{L}_P^*(\cdot)$; finally, we apply these results to the maximum-likelihood estimate $\hat{z}$.

\begin{lemma}\label{lem:refConst}
Consider a $K$-class stochastic blockmodel with membership vector $\bar{z}$, and let $\Pi^z$ denote the partition of its associated $\{P_{ij}\}_{1\leq i<j\leq N}$ induced by any $z \in \{1,\ldots,K\}^N$.  For every $\Pi^z$, there exists a partition $\Pi^*$ that refines $\Pi^z$ and with the property that, if conditions \emph{(i)} and \emph{(ii)} of Theorem~\ref{thm:classConv} hold,
\begin{equation}\label{eq:NerrBndPi}
\bar{L}_P(\bar{z}) - \bar{L}_P^*(\Pi^*) = \frac{N_{\mathrm{e}}(\hat{z})}{N} \, \Omega (M),
\end{equation}
where $N_{\mathrm{e}}(z)$ counts the number of nodes whose true class assignments under $\bar{z}$ are not in the majority within their respective class assignments under $z$.
\end{lemma}

\begin{lemma}\label{lem:refKLDivBnd}
Let $\Pi'$ be a refinement of any partition $\Pi$ of the set $\{P_{ij}\}_{i<j}$; then $\bar{L}_P^*(\Pi') \geq \bar{L}_P^*(\Pi)$.
\end{lemma}

Since Lemma~\ref{lem:refConst} applies to any admissible blockmodel assignment vector $z$, it also applies to the maximum-likelihood estimate $\hat{z}$ for any realization of the data; each $\hat{z}$ in turn induces a partition $\Pi\,\hat{}\,$ of blockmodel edge probabilities $\{P_{ij}\}_{i<j}$, and~\eqref{eq:NerrBndPi} holds with respect to its refinement $\Pi^*$.  By Lemma~\ref{lem:refKLDivBnd}, $\bar{L}_P^*(\Pi\,\hat{}\,) \leq \bar{L}_P^*(\Pi^*)$.  Finally, observe that $\bar{L}_P(\hat{z}) = \bar{L}_P^*(\Pi\,\hat{}\,)$ by the definition of $\bar{L}_P^*$, and so $\bar{L}_P(\bar{z}) - \bar{L}_P(\hat{z}) \geq \bar{L}_P(\bar{z}) - \bar{L}_P^*(\Pi^*)$, thereby establishing~\eqref{eq:NerrBnd}.
\end{proof}

\begin{proof}[of Lemma~\ref{lem:refConst}]
The construction of $\Pi^*$ will take several steps.  For a given membership class under $z$, partition the corresponding set of nodes into subclasses according to the true class assignment $\bar{z}$ of each node.  Then remove one node from each of the two largest subclasses so obtained, and group them together as a pair; continue this pairing process until no more than one nonempty subclass remains, then terminate.  Observe that if we denote pairs by their node indices as $(i,j)$, then by construction $z_i = z_j$ but $\bar{z}_i \neq \bar{z}_j$.

Repeat the above procedure for each class under $z$, and let $C_1$ denote the total number of pairs thus formed. For each of the $C_1$ pairs $(i,j)$, find all other distinct indices $k$ for which the following holds:
\begin{equation} \label{eq:proof3.6}
D\Big(P_{ik} \mid\mid  \frac{P_{ik}+P_{jk}}{2}\Big) + D\Big(P_{jk} \mid\mid  \frac{P_{ik}+P_{jk}}{2}\Big)
\geq C \frac{MK}{N^2},
\end{equation}
where $C$ is the constant from condition (ii) of Theorem~\ref{thm:classConv}, and indices $ik$ and $jk$ in~\eqref{eq:proof3.6} are to be interpreted respectively as $ki$ whenever $k<i$, and $kj$ whenever $k<j$.  Let $C_2$ denote the total number of distinct triples that can be formed in this manner.

We are now ready to construct the partition $\Pi^*$ of the probabilities $\{P_{ij}\}_{1\leq i<j \leq N}$ as follows: For each of the $C_2$ triples $(i,j,k)$, remove $P_{ik}$ (or $P_{ki}$ if $k<i$) and $P_{jk}$ (or $P_{kj}$) from their previous subset assignment under $\Pi^z$, and place them both in a new, distinct two-element subset.  We observe the following:

(i) The partition $\Pi^*$ is a refinement of the partition $\Pi^z$ induced by $z$: Since nodes $i$ and $j$ have the same class label under $z$ in that $z_i = z_j$, it follows that for any $k$, $P_{ik}$ and $P_{jk}$ are in the same subset under $\Pi^z$.

(ii) Since for each class at most one nonempty subclass remains after the pairing process, the number of pairs is at least half the number of misclassifications in that class.  Therefore we conclude $C_1 \geq N_{\mathrm{e}}(z)/2$.

(iii) Condition (ii) of Theorem~\ref{thm:classConv} implies that for every pair of classes $(a,b)$, there exists at least one class $c$ for which~\eqref{eq:proof3.6} holds eventually.  Thus eventually, for any of the $C_1$ pairs $(i,j)$, we obtain a number of triples at least as large as the cardinality of class $c$.  Condition (i) in turn implies that the cardinality of the smallest class grows as $\Omega(N/K)$, and thus we may write $C_2 = C_1 \, \Omega(N/K)$.

We can now express the difference $\bar{L}_P(\bar{z}) - \bar{L}_P^*(\Pi^*)$ as a sum of nonnegative divergences $D(P_{ij} \mid\mid \bar{\theta}_{\zeta_{ij}^*})$, where $\zeta_{ij}^*$ is the assignment mapping associated to $\Pi^*$, and use~\eqref{eq:proof3.6} to lower-bound this difference:
\begin{equation*}
\qquad\qquad\qquad
\bar{L}_P(\bar{z}) - \bar{L}_P^*(\Pi^*) = \sum_{i<j} D(P_{ij} \mid\mid \bar{\theta}_{\zeta_{ij}^*}) = C_2 \, \Omega\Big(\frac{MK}{N^2}\Big) = \frac{N_{\mathrm{e}}(z)}{2} \, \Omega\Big(\frac{M}{N}\Big).
\qquad\qquad\quad\!\QEDlogo
\end{equation*}
\end{proof}

\begin{proof}[of Lemma~\ref{lem:refKLDivBnd}]
Let $\Pi'$ be a refinement of any partition $\Pi$ of the set $\{P_{ij}\}_{i<j}$, and given $a \in \{1,\ldots, L'\}$ indexing $S_a'$, let $F(a)$ denote its index under $\Pi$.  We show that $\bar{L}_P^*(\Pi') \geq \bar{L}_P^*(\Pi)$ as follows:
\begin{align*}
\bar{L}_P^*(\Pi') & = \sum_{a=1}^{L'} |S_a'|\Big\{\bar{\theta}_a' \log \bar{\theta}_a' + (1-\bar{\theta}_a')\log (1-\bar{\theta}_a')\Big\} \\
& \geq \sum_{a=1}^{L'} |S_a'| \Big\{ \bar{\theta}_a' \log \bar{\theta}_{F(a)} + (1-\bar{\theta}_a')\log (1-\bar{\theta}_{F(a)}) \Big\} \\
& = \sum_{b=1}^{L} |S_b| \Big\{ \bar{\theta}_b \log \bar{\theta}_b + (1-\bar{\theta}_b)\log (1-\bar{\theta}_b) \Big\} = \bar{L}_P^*(\Pi),
\end{align*}
where the first inequality holds by nonnegativity of Kullback--Leibler divergence, and the second equality follows from the fact that $\Pi'$ is a refinement of $\Pi$.
\end{proof}

\bibliographystyle{biometrika}
\bibliography{ChoiWolfeAiroldi}

\end{document}